\UseRawInputEncoding

[1994/12/01]
\documentclass[manuscript]{amsart}
\usepackage{amsmath,amsthm,amssymb}

\usepackage{enumitem}

\usepackage{graphicx}

\usepackage{cite}

\newcommand{\supp}{\mbox{supp}}

\makeatletter
\@namedef{subjclassname@2010}{%
  \textup{2010} Mathematics Subject Classification}
\makeatother

\usepackage[T1]{fontenc}

\newtheorem{theorem}{Theorem}[section]
\newtheorem{corollary}[theorem]{Corollary}
\newtheorem{lemma}[theorem]{Lemma}
\newtheorem{proposition}[theorem]{Proposition}
\newtheorem{question}[theorem]{Question}

\theoremstyle{definition}
\newtheorem{definition}[theorem]{Definition}

\numberwithin{equation}{section}

\chardef\bslash=`\\ 

\theoremstyle{definition}

\newcommand{\eval}[2][\right]{\relax
  \ifx#1\right\relax \left.\fi#2#1\rvert}

\title[
On estimate of  operator   for $0<p<\infty $
]
{
 On estimate of  operator   for $0<p<\infty $
 }
\author
[S C Long]
{Shunchao Long}
\address{School of Mathematics and Computational Science\\ Xiangtan University\\ Xiangtan 411105, P.~R.~China.}
\email{sclongc@126.com}

\begin{document}

\begin{abstract}

 Operators such as   Carleson   operator are known to be bounded on  $L^p$ for all $1<p<\infty$, but not from $L^1$ to weak-$L^1$ and from $H^p$ to $L^p$ for each $0<p\leq 1$,  the object of this article is to give a  estimate for all $0<p<\infty$.

 For the  weights $w$  satisfying the doubling condition of order $q$ with $0<q<p$ and the reverse H\"{o}lder condition, by using some new functions spaces,
   we prove that:
\par  $\bullet$ some sublinear operators are bounded from some subspaces of  $L^p_w$ to $L^p_w$ and to themselves  for all $0<p< \infty$; in particular,   these  imply the endpoint estimates
 from  $H^p_w$ to $L^p_w$  and  from  $H^p_w$ to itself for all $0<p\leq 1$;
 these results are applied to many  operators,  such as  Hardy-Littlewood maximal operator, singular integral operators with rough kernels, Calder\'{o}n commutators,  Carleson   operator, the polynomial Carleson operator,
  et al,
 and give the endpoint versions of classical theorems such as Carleson-Hunt theorem and  a conjecture of Stein;

\par  $\bullet$  $H^p_w$ with $0<p\leq 1$  is   characterized by blocks  without vanishing moment conditions;

\par  $\bullet$  $H^p_w$ with $0<p\leq 1$ is   characterized  by a convolution maximal function with a  non-smooth kernel.

 \end{abstract}



\maketitle

\par \section
 { Introduction }

\par \subsection {Endpoint estimate}

The  classical results,
 proved by Hardy and Littlewood  in \cite{HL} and Wiener in \cite{Wie},  Calder\'{o}n and Zygmund in \cite{CZ}, and Carleson  in \cite {Carl} and Hunt in \cite{Hunt}, respectively,
  state that
 Hardy-Littlewood maximal operator and  the singular integral operators  on ${\bf R}^n$,   and  Carleson  operator  on ${\bf R}^1$  are  bounded from  $L^p$ to itself for $1<p<\infty$.
Unfortunately,   these do not hold  for $0<p\leq 1$, (see  page 79 and page 284 of Grafakos \cite{Grafakos}, and Kolmogorov \cite{Ko1}).
As a substitute, for $p=1$, the weak type (1,1) estimate from  $L^1$ to weak-$L^1$ is considered.
 It is proved that Hardy-Littlewood maximal operator and  the singular integral operator are  of weak type (1,1), (see \cite{HL, Wie} and \cite{CZ}), in fact, it is an important step for $L^p$ estimate for $1<p<\infty$.  But
  Carleson  operator is not  of weak type (1,1),
 (see  Kolmogorov \cite{Ko1}). At the same time, as another substitute, for $0<p\leq 1$, the estimate from  $H^p$ to  $L^p$ (or, to  itself)  is considered.
   In \cite{FS72}, Fefferman and Stein  first proved  the $H^p$ boundedness of  some singular integral operators  for some $0<p\leq 1$.    Garc\'{\i}a-Cuerva and  Rubio de Francia proved in \cite{GRubio} that some regular singular integral operators of principle value type (including Hilbert transform) are bounded from $H^p$ to  $L^p$ and  from $H^p$ to  itself  for  $n/(n+1)<p\leq 1$.
  However,  for  $0<p \leq n/(n+1)$, this does not hold, (see \cite{Grafakos}).  While, Hardy-Littlewood maximal operator and  Carleson   operator  fail to be bounded from  $H^p$ to  $L^p$ for all $0<p\leq 1$  (see   \cite{Grafakos} and Stein \cite{Stei}). Carleson   operator  is not even bounded from  $H^1$ to  weak-$L^1$,
  (see Zygmund  \cite{Zy},  Chapter 8).

   At the same time,  the weighted version of  the results mentioned above   holds for the weights in  the classical  Muckenhopt  class  $A_p$.
    Muckenhopt   \cite{Muck},  Hunt,  Muckenhoupt and  Wheeden   \cite{HuMu},  Coifman and  Fefferman   \cite{CoFe}, and    Hunt and Young  \cite{HY} proved,  respectively,  that  Hardy-Littlewood maximal operator, the singular integral operators,   and  Carleson  operator  are  bounded on $L_w^p$   for $1<p<\infty$ and $w\in A_p$,  Hardy-Littlewood maximal operator and the singular integral operators are also  bounded from $L_w^1$ to weak-$L_w^1$ for   $w\in A_1$, (see  \cite{Muck} and \cite{CoFe}). The regular singular integral operators of principle value type are bounded from $H_w^p$ to  $L_w^p$ and  from $H_w^p$ to  itself  for  $n/(n+1)<p\leq 1$ and $w\in A_1$, (see \cite{LL}).
    And the negative results under unweighted case are still true for weighted case since   $1\in A_1$, i.e., there is $w\in A_1$ such that,  above operators    fail to be bounded from  $L_w^1$ to  itself,
    Hardy-Littlewood maximal operator and  Carleson  operator  fail to be bounded from  $H_w^p$ to  $L_w^p$ for all $0<p\leq 1$,
   and the regular singular integral operators of principle value type fail to be bounded from  $H_w^p$ to  $L_w^p$ for   $0<p \leq n/(n+1)$.

These cases appear in the estimates of lots of operators in harmonic analysis.
 Many  operators   have been proved to be bounded on  $L^p$  (or $L_w^p$  for $w\in A_p$) for $1<p<\infty$,
 (see, e.g. \cite{St61970,S3,GRubio,Mey,Grafakos,HeoNS ,Lie,F2,Chanillo,Chen,DR,KrL,PhongStein,RicciStein,SS,CZ2,Cald,ChristJ,Muckenhoupt2,BCo,HeoHY}).
  Generally, these operators are not bounded on  $L^1$  (or $L_w^1$  for $w\in A_1$).  A lot of  efforts were devoted to the endpoint estimate for $0<p\leq 1$: weak type (1,1) estimate
  (see e.g.  \cite{St61970,S3,GRubio,Mey,Grafakos,ChristF, Hofmann, Seeger, Seeger2, Tao3, DingL} ) and $H^p$ estimates
  (see e.g.  \cite{S3,GRubio,Mey,Grafakos,F2, Christ1, HeoNS,STW,STW2,SeegerW,SeegerT,Tao,Tao2,HeoHY} ).
 For these operators, some  are proved to be of weak type (1,1) estimate  but fail to be of $H^p$ estimate  for all $0<p\leq 1$; some are   of $H^p$ estimate  for some $0<p\leq 1$ but fail to be of weak type (1,1); some are    of both weak type (1,1) estimate and $H^p$ estimate  for some $0<p\leq 1$ but  the $H^p$ estimate is   not true for $0<p$ enough small;
  even some are known to be of neither  weak type (1,1) estimate nor $H^p$ estimate for all $0<p\leq 1$;
  while for  those remaining operators, it is open problem whether there are weak type (1,1) estimate and $H^p$ estimate for all $0<p\leq 1$,
for example, very recently,  Lie proved in  \cite{Lie} the one dimensional case of a conjecture  of Stein  which state that  the polynomial Carleson operator  is bounded on $L^p$ for any $1<p<\infty$ (see  \cite{Stein4} and  \cite{Stein5}),
 while the problems of  endpoint estimate for $0<p\leq 1$  (including  both $H^p$ estimate  and  weak type (1,1) estimate) are open. In addition, it is also worth noting that the maximal operators are obviously   not bounded from $H ^ p $ to itself for all $0<p\leq 1$  since $H^p$ functions have the vanishing properties   \cite{S3}.

In short,   for a lot of operators, for $0<p$ enough small, we don't seem to have any estimates  similar to their $L^p$ estimates with $1<p<\infty$.

   It would be natural to ask the following question.
\begin{question}
For a operator  mentioned above, is   there  a    endpoint estimate
   for $0<p$ enough small? further more, is   there  a  estimate for all $0<p<\infty$  similar to  $L^p$ estimate with $1<p<\infty$?
\end{question}
\par The first purpose of the present paper is to give an answer to  this question.

 For all $0<p\leq 1 $, we will give $H^p_w$ estimates for some operators, these imply the weighted endpoint versions  of the famous theorems  of Hardy, Littlewood and Wiener  in \cite{HL, Wie} (for Hardy-Littlewood maximal operator), of Calder\'{o}n and Zygmund in \cite{CZ2} (for  the singular integral operators with rough kernel), of Calder\'{o}n in \cite{Cald} (for Calder\'{o}n commutator), of  Carleson and Hunt  in \cite {Carl,Hunt} (for  Carleson  operator), and  a conjecture of Stein in  \cite{Stein4,Stein5} (for  the polynomial Carleson operator), et.al.. In fact, we  will  give  a   new weighted estimate for all $0<p<\infty$ for these operators. (See Section 2).

  \par \subsection {Characterization of Hardy space}
  \par \subsubsection { Maximal functions characterization}
  $H^p$, $p>0$, consists of those tempered distributions $f \in {\mathcal{S}}'$ for
which the maximal function $M_{\varphi}f(x)=\sup_{t > 0}|\varphi_t*f(x)|\in L^p$,
 where  $\varphi $ is a function in ${\mathcal{S}}$, the Schwartz space of rapidly decreasing smooth functions,
 satisfying $\int_{{\bf R}^n}\varphi (x) dx=1$,
 and $\varphi_t(x)=t^{-n}\varphi(x/t),  t>0, x\in {\bf R}^n$, (see  Fefferman and Stein \cite{FS72}, see also \cite{S3}). In \cite{W79}, Weiss considered the problem reducing  smoothness of $\varphi$
   in Hardy space theory.  Unfortunately,  the example $\tilde{\varphi}=\frac{1}{|B(0,1)|}\chi_{B(0,1)}$   shows that the assumption of smoothness of  $\varphi$  cannot be removed in the characterization of $H^1$,  since $\{f\in L^1:M_{\tilde{\varphi}}f\in L^1\}=\{0\}$.

 \par \subsubsection { Atoms characterization}
At the same time, Coifman in  \cite{Co1} and Latter in  \cite{ Lat}  proved that $H^p, 0<p\leq 1$, can  also be characterized   in terms of
atoms   satisfying
   a compact support condition,
  a size  condition and
   some vanishing moment conditions.
Attempts to reduce  vanishing moment conditions in the atom theory  of Hardy space appear in many literatures,
  (see, for example, Stein \cite{S3}).
   However, it is known that a bounded, compactly supported function $f$ belongs to $H^p$ if and only if it satisfies the vanishing  moment conditions $\int x^{\alpha}f(x)dx=0$ for all $|\alpha|\leq n(p^{-1}-1),$ (see \cite{S3},p129). This shows that the  vanishing moment conditions of $H^p$ atoms cannot be removed in the classical $H^p$ theory.

\par

\par \subsubsection {Characterization of weighted Hardy space}
For a weight $w$ in  the classical  Muckenhopt  class, the weighted Hardy space    $H_w^p$ with  $0<p\leq 1$ is characterized by the convolution  maximal functions with  smooth kernels and  by atoms  with vanishing moment conditions, (see Garc\'{\i}a-Cuerva  \cite{Gc} and Str\"{o}mberg and Torchinsky  \cite{ST}). But, since $1\in A_1$, we see from Subsection 1.2.1 and 1.2.2  that   there is $w\in A_1$ such that   $H_w^p$ with $0<p\leq 1$ can neither be characterized by convolution  maximal function with a non-smooth kernel nor blocks  without vanishing moment conditions.

Naturally, we ask the following questions:
\begin{question}
  Is   there  a class  of  weights such that, for all weights $w$ in this class,
   $H^p_w$  with $0<p\leq 1$ is    characterized by   convolution maximal function with a  non-smooth kernel?
\end{question}

\begin{question}
  Is   there  a class  of  weights such that, for all weights $w$ in this class,
 $H^p_w$ with $0<p\leq 1$ is   characterized by blocks  without vanishing moment conditions?
\end{question}

\par We will also  give
affirmative
answers to  these two questions.

\par \subsection {Main theorems}

 We write  $\bar{p}=\inf \{p,1\}$ for $0<p<\infty$.  A  nonnegative local integrable function is called a weight.

  Our solution to the above three questions is based on the following new functions spaces generated by blocks.

\par
\begin{definition}\label{def_3.1}
Let $  0< p < \infty ,
0< s \leq \infty,  w$ be a weight.

 \par A function $a$ is said to be an
 $ (p, s, w)$-block, if there is a cube $Q$ in ${\bf R}^n$ such that

\par (i)~~~~supp $ a\subseteq Q,$
\par (ii)~~~~$\|a\|_{L^{s}   }\leq |Q|^{1/s}w(Q)^{-1/p}$.

\par The spaces $B^{p,s}_{w}$ consists of functions $f$ which can be written  as
$f=\sum\limits_{k=1}^{\infty} \lambda _k h_k ,$
   where   $h_k$ are $(p,s, w)$-blocks and $\lambda _k$ are real numbers with
 $\sum\limits_{k=1}^{\infty} |\lambda _k|^{\bar{p}} <+ \infty .$

 \par We equip  $B^{p,s}_{w}$ with the  quasi-norm
$
\|f\|_{B^{p,s}_{w}}=
 \inf  \left( \sum\limits_{k=1}^{\infty} |\lambda _k|^{\bar{p}} \right)^{1/\bar{p}},
  $
where the infimum is taken over all the above decompositions of $f$.
\end{definition}

For the power weight  $w(x)=|x|^{\alpha}$ with $-n<\alpha<n(p-1)$, these spaces were introduced by the author in \cite{Longs,Longs2}. The unweighted block spaces with different quasi norm were introduced by M.Taibleson and G.Weiss to study a.e. convergence of Fourier series   in \cite{TW1}, see also  Lu, Taibleson and Weiss  \cite{LTW1}.

We need also the following definitions.
\begin{definition}
  Let   $0< p<\infty$. We say  a weight $w\in D_{p}$  if    $w$  satisfies    the doubling condition of order $p$
 \begin{equation}\label{2.2}
 w(\lambda Q)\leq C  \lambda^{np} w(Q)
  \end{equation}
   for any cube $Q $ and $\lambda>1$, where $C$ is a  constant  independent of $Q $ and $\lambda$.
\end{definition}

\par
\begin{definition}
  Let   $1< r<\infty$. We say a weight
 $w \in RH_r$  if $w$
 satisfies   the reverse H\"{o}lder condition of order $r$
 \begin{equation}\label{2.2}
 \left(\frac{1}{|Q|}\int_Q w^{r}\right)^{1/r}\leq \frac{C}{|Q|}\int_Q w
\end{equation}
 for  every cube $Q$,  where  $C$ is a  constant  independent of $Q$.

\end{definition}

\begin{definition}\label{def_2.1}
 We say a weight $w\in  {P} $ if there exist a sequence $\{Q_i\}$ of cubes with $w(Q_i)>0$ for each $i$ and  whose interiors are disjoint each other such that  ${\bf R}^n=\bigcup_{i=1}^{\infty} Q_i$.

\end{definition}

Recall that
the classical weighted  Hardy spaces $H_w^p$, $p>0$, $w\in A_{\infty}$, consists of those tempered distributions $f \in {\mathcal{S}}'$ for
which the maximal function
$$M_{\varphi}f=\sup_{t > 0}|\varphi_t*f(x)|\in L_w^p,$$
 for some   $\varphi $  in ${\mathcal{S}}$,
 satisfying $\int_{{\bf R}^n}\varphi (x) dx=1$,
 and $
\varphi_t(x)=t^{-n}\varphi(x/t),  t>0, x\in {\bf R}^n
$,  (see \cite{Gc,ST}).

\par
 We denote $ B^{p,s}_{w}$ by $ BL^{p,s}_{w}$,  when $ 0< p <  \infty$,
   $w\in RH_r$  with $ 1<r <\infty  $,
and  $   rp/(r-1)\leq s \leq \infty $, or,
 $w$ be a weight and  $  s=\infty$, (since $ B^{p,s}_{w}\subset L^{p}_{w}$, see Proposition 3.3).

\par  We denote $ B^{p,s}_{w}$ by $ BH^{p,s}_{w}$, when $ 0< p <  \infty$, $w\in A _{q,r} $ with $ 0<q <p  $  and  $ 1<r <\infty  $,
  and $   1 < s\leq \infty  $ and $  rp/(r-1) \leq s   $, (since $ B^{p,s}_{w}\subset H^{p}_{w}$, see Theorem 4.1).

\par We denote $D _{q} \cap RH_r $ by  $A _{q,r} $ for $ 0<q <\infty  $  and  $ 1<r <\infty  $.

Our answer to Question 1.3 is as follows.

\begin{theorem}
  Let $ 0< p \leq 1, w\in A _{q,r}     $ with $ 0<q <p  $  and  $ 1<r <\infty  $,
then,
\begin{equation}
H^{p}_{w} = BH^{p,s}_{w}
\end{equation}
 for    $   1 < s\leq \infty  $ and $  rp/(r-1) \leq s   $. It is sharp in the sense that  (1.3) may not hold when $ q= p =1$.

\end{theorem}

\par  Theorem 1.8 gives a block characterization of weighted Hardy spaces $H^{p}_{w}$.

\par
  For a function $f$, we define formally his Hardy-Littlewood maximal function as $Mf(x)=\sup_{x\in Q}\frac{1}{|Q|}\int_Q|f(y)|dy$ for all cubes $Q$.
         Let $w$ be a weight, set, for $0<p\leq \infty$,
 $$ML^{p}_{w}=\{f:Mf\in L_w^p\}.$$

\par Our answer to Question 1.2 is as follows.

\begin{theorem}\label{th_4.7}
Let $0<p \leq 1,    w\in  A _{q,r}      \cap  P $ with $ 0<q <p $  and  $ 1<r <\infty  $,
   then
    \begin{equation}
    H^{p}_{w} = ML^p_{w}.
    \end{equation}
  It is sharp   in the sense that  (1.4) may not hold when $ q= p =1$.
\end{theorem}

\par  Let $\tilde{\varphi}=\frac{1}{|B(0,1)|}\chi_{B(0,1)}, \tilde{\varphi}_t(x)=t^{-n}\tilde{\varphi}(x/t),  t>0,$
    and  define   $M_{\tilde{\varphi}}f=\sup_{t > 0}|\tilde{\varphi}_t*f(x)|$ formally for  a function $f$. Set
 $M_{\tilde{\varphi}}L^{p}_{w}=\{f:M_{\tilde{\varphi}}f\in L_w^p\}.$ It is easy to get $M_{\tilde{\varphi}}L^{p}_{w}=ML^{p}_{w}$, then,
Theorem 1.9 gives a characterization of $H^{p}_{w}$
by  a convolution maximal function with   non-smooth kernel,
 since $\tilde{\varphi}$ is not a smooth function.

\par For Question 1.1, we will prove for $ w\in  A _{q,r}$ with $0<q<p$ and $1<r<\infty$ that

   $ many~ operators~      is~ bounded~ from $ $H_w^p$ $to$  $L_w^p$ $and~ from$ $H_w^p$ $to~  itself~  for ~all$ $0<p\leq 1$.

  These are the special cases of the following results:

 $ these~ operators~ is~ bounded~ from ~  some ~subspaces ~of~L_w^p~to~L_w^p~and~ from~$ $ some ~ subspaces~ of~L_w^p~to  ~itself~ for ~all~0<p<\infty$.

In fact, we will prove the above  results for the following  sublinear operator $T$ which includes the above  operators.

   $T$ is defined   for every $(p,s,w)-$bolck $h$ with supp $h \subset Q$, a cube with the cental  $x_0$, and  satisfies  the size
condition
\begin{equation}\label{7.1}
|Th(x)| \leq  C \frac{  \|h\|_{L^{1}}}{|x-x_0|^n},
\end{equation}
for almost all $x\in (2n^{1/2}Q)^c $,
and
\begin{equation}\label{7.3}
 \|Th\|_{L^s(2n^{1/2}Q)} \leq C \|h\|_{L^s(Q)},
\end{equation}
for some $1\leq s\leq \infty$.
  $T$  is either a linear operator or a sublinear operator satisfying the following condition
\begin{equation}\label{7.2}
|Tf(x)|\leq \sum|\lambda_j||Ta_j(x)|,~w{\rm -a.e.},
\end{equation}
 for every $f=\sum\lambda_ja_j$ in $  {B}L^{p,s}_{w}$.

 \par We say a proposition holds $w$-a.e. if there  is a set $E$ with $w(E)=0$ such that this proposition holds for all $x\in E^c$.

\par We say $f=\sum\lambda_ja_j$ in $  {B}L^{p,s}_{w}$   means that    each  $a_j$ is a $(  p,s, w)$-block and $\sum|\lambda_j|^{\bar{p}}<\infty$.

\par The conditions (1.5), (1.6) and (1.7) are satisfied by many   operators in harmonic analysis (see Subsection 7.3).

\par Our answer to the question 1.1 is as follows.

\begin{theorem}\label{th_7.1}
    Let $0<p <\infty $. Let
    \par  (a) $ w\in D _{q}     $ with $ 0<q <p  $,  and $ s = \infty $, or
   \par (b)
  $ w\in A _{q,r}     $ with $ 0<q <p  $  and  $ 1<r <\infty  $,  and $ \max \{ rp/(r-1),1\}\leq s \leq \infty$.
     \\Suppose that an operator $ T $ is defined for every $(p,s,w)$-block and satisfies (1.5)  and (1.6).

\par     (i)  If   $ T $  is a linear operator, then
  $T$ has an unique bounded extension (still denoted by $T$) from $BL^{p,s}_{w}$
    to $L^{p}_{w}$ that satisfies
 \begin{equation}\label{7.4}
  T(\sum_{j=1}^{\infty}\lambda_ja_j)(x) =\sum_{j=1}^{\infty}\lambda_jTa_j(x),
  \end{equation}
  in $L^p_w$ and $w$-a.e., and
  \begin{equation}\label{7.4}
  \|Tf\|_{L^{p}_{w}}   \leq   C \|f\|_{BL^{p,s}_{w}},
    \end{equation}
    for all  $ f= \sum_{j=1}^{\infty}\lambda_ja_j \in BL^{p,s}_{w}$.

   \par (ii)  If  $ T $ is a sublinear operator
satisfying (1.7)  for    $f=\sum\lambda_ja_j \in BL^{p,s}_{w} $,
 then $T$ is bounded from $BL^{p,s}_{w}$
   to $L^{p}_{w}$  and (1.9)
  holds.

\end{theorem}

\par In particular, by Theorem 1.10 and Theorem 1.8, we have
\begin{corollary}\label{th_7.1}
    Let $0<p \leq 1, w\in A _{q,r}     $ with $ 0<q <p  $  and  $ 1<r <\infty  $.  Suppose that an operator $ T $ is defined for every $(p,s,w)$-block with   $   1 < s\leq \infty  $ and $  rp/(r-1) < s   $
    and satisfies (1.5)  and (1.6).

\par     (i)  If  $ T $ is a linear operator, then
  $T$ has an unique bounded extension (still denoted by $T$) from $H^{p}_{w}$
    to $L^{p}_{w}$ that satisfies (1.8)
   in $H^p_w$ and $w$-a.e., and
  \begin{equation}\label{7.4}
  \|Tf\|_{L^{p}_{w}}   \leq   C \|f\|_{H^{p}_{w}},
    \end{equation}
    for all  $ f= \sum_{j=1}^{\infty}\lambda_ja_j \in H^{p}_{w}=BL^{p,s}_{w}$.

   \par (ii)  If  $ T $
 is a sublinear operator
satisfying (1.7)  for    $f=\sum\lambda_ja_j \in H^{p}_{w} $,
 then $T$ is bounded from $H^{p}_{w}$
   to $L^{p}_{w}$  and (1.10)
  holds.

\end{corollary}

\begin{theorem}\label{th 7.6}
       Let $0<p <\infty, w\in A _{q,r}     $ with $ 0<q <p  $  and  $ 1<r <\infty  $.  Let   $\max\{ rp/(r-1), p/q\}< s \leq \infty $.
 Suppose that an operator $ T $ is defined for every $(p,s,w)$-block,  satisfies (1.5)  and is bounded on $L^s$.

\par     (i)  If   $ T $  is a linear operator, then
  $T$ has an unique bounded extension   (still denoted by $T$)  from $BH^{p,s}_{w}$ to itself
   that satisfies (1.8) in $BH^{p,s}_{w}$, and
  \begin{equation}\label{7.5}
  \|Tf\|_{BH^{p,s}_{w}}   \leq   C \|f\|_{BH^{p,s}_{w}}.
   \end{equation}
  for all  $ f= \sum_{j=1}^{\infty}\lambda_ja_j \in BH^{p,s}_{w}$.

\par (ii)  If $ T $ is a sublinear operator
satisfying (1.7)  for    $f=\sum\lambda_ja_j \in BH^{p,s}_{w} $,
 and  $   w\in    P   $,
 then $T$ is bounded from $  BH^{p,s}_{w}$ to itself
   and (1.11)   holds.

\end{theorem}
In particular, by Theorem 1.12 and Theorem 1.8, we have
\begin{corollary}\label{th 7.6}
       Let $0<p \leq 1, w\in A _{q,r}     $ with $ 0<q <p  $  and  $ 1<r <\infty  $.
 Suppose that an operator $ T $ is defined for every $(p,s,w)$-block with $\max\{ rp/(r-1), p/q\}< s \leq \infty $,
  satisfies (1.5)    and is bounded on $L^s$.

\par     (i)  If  $ T $  is a linear operator, then
  $T$ has an unique bounded extension   (still denoted by $T$)  from $H^{p}_{w}$ to itself
   that satisfies (1.8) in $H^{p}_{w}$, and
  \begin{equation}\label{7.5}
  \|Tf\|_{H^{p}_{w}}   \leq   C \|f\|_{H^{p}_{w}}.
   \end{equation}
  for all  $ f= \sum_{j=1}^{\infty}\lambda_ja_j \in H^{p}_{w}=BH^{p,s}_{w}$.

\par (ii)  If $ T $ is a sublinear  operator satisfying (1.7)  for    $f=\sum\lambda_ja_j \in H^{p}_{w} $, and $ w\in   P   $,
 then $T $ is bounded from $H^{p}_{w}$ to itself    and (1.12)   holds.

\end{corollary}

 Theorem 1.10 and Theorem 1.12 hold for $w\in A_q$ with $1\leq q < p\leq \infty$, since there is $1<r<\infty$ such that $A_q\subset  A _{q,r}     $ for $1\leq q <\infty$.

\section{Applications }

Theorem 1.10 and Theorem 1.12 are applied to many   operators in harmonic analysis.

Let  $T$ satisfy the following size condition:
\begin{equation}\label{7.11}
|Tf(x)|\leq C \int_{{\bf R}^n}\frac{|f(y)|}{|x-y|^n}dy, ~~~~ x\notin
{\rm supp}f,
\end{equation}
for any  integral function $f$ with compact support.
 We have that (2.1) implies (1.5), (see Lemma 7.1). At the same time, the $L^s$ boundedness of $T$ obviously implies (1.6). Therefore, we have
\begin{theorem}\label{}
(i) Theorem 1.10   holds if the conditions (1.5) and (1.6) are replaced by (2.1) and the $L^s$ boundedness of $T$, respectively.
\par (ii) Theorem 1.12 holds if the conditions (1.5) is replaced by (2.1).
\end{theorem}
(2.1) is satisfied by many operators  (see \cite{SW}). Thus, Theorem 1.10 and Theorem 1.12 are applied to these   operators.

\subsection{Applications to linear operators }
 Theorem 1.10(i) and Theorem 1.12(i) are applied to the following   linear operators.

\par  $\bullet$  {\bf Hilbert transform}
 $$Hf(x)={\rm p.v.} \int_{\bf R}\frac{f(y)}{x-y}dy
 . $$

 $H$ is bounded from  $L^p_w$  to itself for all $1<p<\infty$ and $w\in A_p$  and from $L^1_w$  to weak-$L^1_w$  for   $w\in A_1$, ( see Hunt, Muckenhoupt and  Wheeden  \cite{HuMu}, and see Stein\cite{St61970} for the case $w=1$ ), and   from  $H^p_w$ to itself for all $1/2<p\leq 1$ and $w\in A_1$,  ( see Lee and Lin \cite{LL} ). But, $H$ fails to be  bounded both from  $L^1$  to itself and from  $H^p$  to   $L^p$ for $0<p\leq 1/2$, (see \cite{Grafakos}).

\par $\bullet$  {\bf Riesz transform}
$$R_j f(x)={\rm p.v.} \int_{{\bf R}^n}\frac{|x_j-y_j|}{|x-y|^{n+1}} f(y)dy, j=1,2,\cdots,n.$$

   $R_j$ is bounded from $L^p_w$ to itself for all $1<p<\infty$ and $w\in A_p$  and from $L^1_w$  to weak-$L^1_w$  for   $w\in A_1$, (see Coifman and Fefferman \cite{CoFe}, and see Stein \cite{St61970} for the case $w=1$ ),   and from  $H^p_w$ to itself for all  $n/(n+1)<p\leq 1$ and $w\in A_1$,   ( see Lee and Lin \cite{LL} ). But $R_j$ fails to be  bounded both from  $L^1$  to itself and from  $H^p$  to   $L^p$ for $0<p\leq n/(n+1)$, (see \cite{Grafakos}).

\par $\bullet$  {\bf Calder\'{o}n-Zygmund  operator $T_{CZ} $}.

 Let $\Delta =\{(x,x):x\in {\bf R}^n $ be the diagonal of ${\bf R}^n \times {\bf R}^n$. We define a  Calder\'{o}n-Zygmund kernel $K(x,y)$ to be a locally integrable function ${\bf R}^n \times {\bf R}^n\setminus\Delta\rightarrow {\bf R}$ which satisfies
\begin{equation} |K(x, y)| \leq C\frac{1}{|x - y|^{n} },\end{equation}
and
\begin{equation}|K(x, y)-K(x, z)| +|k(y, x)-k(z, x)| \leq C\frac{|y-z|^{\delta}}{|x - z|^{n+\delta} }\end{equation}
for  $2 |y - z| < |x - y|$ and some $ 0<\delta \leq 1$.
Let $T_{CZ}:C^{\infty}_0\rightarrow \mathcal{D}'$ be a bounded linear operator.
 $T_{CZ}$ is called a Calder\'{o}n-Zygmund operator if $T_{CZ}$ extends to  a  bounded operator on $L^2$ and there exist a  Calder\'{o}n-Zygmund kernel $k(x,y)$ such that
\begin{equation*}T_{CZ}^{}f(x) = \int K(x,y)f(y)dy\end{equation*}
 for any $f\in C^{\infty}_0$ and $x\in \supp (f)$.

 $T_{CZ}$ is bounded from $L^p_w$ to itself for all $1<p<\infty$ and $w\in A_p$ and from $L^1_w$ to weak-$L^1_w$ for $w\in A_1$.
But $T_{CZ}$   fails to be bounded on  $L^1$.
$T_{CZ}$   is bounded  from  $H^p$ to  $L^p$  for $n/(n+\delta)<p\leq 1$,
but,   it may not be bounded from   $H^p$  into itself for $n/(n+\delta)<p\leq 1$, since
$T_{CZ}$ is bounded from   $H^p$  into itself
      if and only if $\widetilde{T}1=0$, where $\widetilde{T}$ is the conjugate operator of $T$,
At the same time, the boundedness of $T_{CZ}$ from  $H^p$ to  $L^p$  may not be true for $0<p\leq n/(n+\delta)$. (See Meyer \cite{Mey}).

\par $\bullet$   {\bf Singular integral operators with rough kernel $T_{\Omega}$}.

Let  $ \Omega $ be a function defined on ${\bf R}^n \backslash {0}, n\geq 2, $  satisfying
    \begin{equation}\label{}
     \Omega (rx')= \Omega (x') ~for~any~r>0~ and ~x'\in S^{n-1},~and~
     \Omega \in L^\infty (S^{n-1}).
    \end{equation}
    Define
    $$T_{\Omega} f(x)={\rm p.v.} \int_{{\bf R}^n}\frac{\Omega (x-y)h(|x-y|)}{|x-y|^{n}} f(y)dy, $$
    where
  \par (a) $h=1$, and  $ \Omega $ is odd and satisfies (2,4); or,
  \par (b) $ h\in L^\infty ([0,\infty))$, and  $ \Omega $  satisfies (2,4) and
   \begin{equation}\label{}
       \int _{{\bf S}^{n-1}}\Omega(\theta)d\theta=0,
     \end{equation}
 where $d\theta$ denotes the surface measure of ${\bf S}^{n-1} $.

For the case (a),  by Calder\'{o}n and Zygmund   \cite{CZ2}, $T_{\Omega}$ is bounded from $L^p$ to itself for all $1<p<\infty$;
 for the case (b), $T_{\Omega}$ is bounded from $L^p_w$ to itself for all $1<p<\infty$ and $w\in A_p$, (see Duoandikoetxea and Rubio de Francia \cite{DR}, and see Chen \cite{Chen} for the case $w=1$).
 For $ h=1$, if $ \Omega $ satisfies (2.4) and (2.5),
 $T_{\Omega}$ is bounded from $L^1$  to weak-$L^1$, (see Christ  \cite{ChristF}, Hofmann \cite{Hofmann} and Seeger \cite{Seeger}).

\par $\bullet$  {\bf General singular integral of Muckenhoupt type}
 $$T_{\Omega,ir} f(x)={\rm p.v.} \int_{{\bf R}^n}\frac{\Omega (x-y) }{|x-y|^{n+ir}}  f(y)dy, $$
where   $i=\sqrt{-1} , r\in {\bf R}\backslash \{0\}, \Omega $   satisfies (2,4) and (2.5).

 By B. Muckenhoupt \cite{Muckenhoupt2}, $T_{\Omega,ir}$ is bounded from $L^p$ to itself for all $1<p<\infty$.

\par $\bullet$  {\bf  Calder\'{o}n commutator }
 $$T_{\Omega,A} f(x)={\rm p.v.} \int_{{\bf R}^n}\frac{\Omega (x-y) }{|x-y|^{n}} \frac{A(x)-A(y)}{|x-y|}  f(y)dy, $$
  where $ A\in Lip({\bf R}^n), \Omega $  satisfies (2,4) and one of the following conditions:
  \par (a)  $ \Omega $ is even ;
  \par (b)     $ \Omega $ is add and  satisfies
\begin{equation}\label{}
       \int _{{\bf S}^{n-1}}\Omega(\theta)\theta^{\alpha}d\theta=0, ~~for ~~all~~\alpha\in {\bf Z}_+^n~~with~~|\alpha|=1.
     \end{equation}

By Calder\'{o}n   \cite{Cald},  $T_{\Omega,A}$ is bounded from $L^p$ to itself for all $1<p<\infty$.

\par $\bullet$  {\bf Higher order Calder\'{o}n commutator }
 $$T^k_{\Omega,A} f(x)={\rm p.v.} \int_{{\bf R}^n}\frac{\Omega (x-y) }{|x-y|^{n}}\left(\frac{A(x)-A(y)}{|x-y|}\right)^k f(y)dy, $$
  where $k\geq 1, A\in Lip({\bf R}^n), \Omega $  satisfies (2,4) and
\begin{equation*}\label{}
       \int _{{\bf S}^{n-1}}\Omega(\theta)\theta^{\alpha}d\theta=0, ~~for ~~all~~\alpha\in {\bf Z}_+^n~~with~~|\alpha|=k.
     \end{equation*}

S. Hofmann \cite{Hofmann3} proved that $T^k_{\Omega,A}$ is bounded from $L^p_w$ to itself for all $1<p<\infty$ and $w\in A_p$.

\par $\bullet$  {\bf General Calder\'{o}n commutator}
 $$T_{\Omega,F,A} f(x)={\rm p.v.} \int_{{\bf R}^n}\frac{\Omega (x-y) }{|x-y|^{n}}F\left(\frac{A(x)-A(y)}{|x-y|}\right) f(y)dy, $$
under  the following conditions:
  \par (a)  $ \Omega $ is odd and satisfies (2,4);
 \par (b)  $A\in Lip({\bf R}^n)$;
 \par (c) $ F(t) $ is odd for $t\in {\bf R}$  and  is real analytic in $\{ |t|\leq \|\nabla A\|_{\infty}\}$.

By A. P. Calderon, C. P. Calderon, E. Fabes, M. Jodeit, and N. M. Rivi\`{e}re \cite{CCFJR}, $T_{\Omega,A,l}$ is bounded from $L^p$ to itself for all $1<p<\infty$.

\par $\bullet$  {\bf   Calder\'{o}n commutator  of Baj\u{s}anski-Coifman type}.

Let $A_{\alpha}(x)=\partial_x^{\alpha}A(x), \alpha\in {\bf Z}_+^n$ and
$$
P_l(A,x,y)=A(x)-\sum_{|\alpha|<l}\frac{A_{\alpha}(y)}{\alpha!}(x-y)^{\alpha},
$$
where $l\in {\bf N}$. Define
 $$T_{\Omega,A,l} f(x)={\rm p.v.} \int_{{\bf R}^n}\frac{\Omega (x-y) }{|x-y|^{n}}\frac{P_l(A,x,y)}{|x-y|^l} f(y)dy, $$
where,
   (a)  $ \Omega $   satisfies (2,4) and (2.6);
  (b)  $A_{\alpha} \in L^{\infty}({\bf R}^n)$ for $|\alpha|=l$.

Clearly, when $l=1$, $T_{\Omega,A,1}$ is just the Calder\'{o}n commutator $T_{\Omega,A}$.

By Baj\u{s}anski and Coifman \cite {BCo}, $T_{\Omega,A,l}$ is bounded from $L^p$ to itself for all $1<p<\infty$.

By Ding and Lai \cite {DingL},  $T_{\Omega,ir}, T_{\Omega,A}, T^k_{\Omega,A},T_{\Omega,A,l},T_{\Omega,ir}$ are bounded from $L^1$ to weak-$L^1$.

\par $\bullet$  {\bf Commutator of Christ-Journ\'{e}  type}

Let  $k\in {\bf N}, a\in L^{\infty}$ and $m_{x,y}a=\int_0^1a(sx+(1-s)y)ds$.
Define
\begin{equation*}T_{a,k}^{}f(x) = {\rm p.v.}\int K(x-y)(m_{x,y}a)^kf(y)dy,\end{equation*}
where $K$ be the
Calder\'{o}n-Zygmund convolution kernel.
Define
\begin{equation*}T_{ \Omega ,a,k}^{}f(x) = {\rm p.v.}\int\frac{\Omega (x-y) }{|x-y|^{n}}(m_{x,y}a)^kf(y)dy,\end{equation*}
where  $ \Omega $   satisfies (2,4).

When $n=1$, $m_{x,y}a=\frac{\int_0^x a(z)dz-\int_0^y a(z)dz}{x-y}$. Let $K(x)=\frac{1}{x}$ and $A(x)=\int_0^x a(z)dz$, then, $A'(x)=  a(x)\in L^{\infty}({\bf R})$. Then,
\begin{equation*}T_{a,1}^{}f(x) = {\rm p.v.}\int \frac{A(x)-A(y)}{x-y}\frac{f(y)}{x-y}dy,\end{equation*}
which is the  Calder\'{o}n commutator $T_{\Omega,A}$ with $\Omega=1$.

 M. Christ and J.-L. Journ\'{e} \cite{ChristJ}  proved that $T_{a,k}$ is bounded from $L^p$ to itself for all $1<p<\infty$.  Ding and Lai proved \cite{DingL2} that $T_{a,1}$ is bounded from $L^p_w$ to itself for all $1<p<\infty$ and $w\in A_p$.
Grafakos and P. Honz\'{i}k \cite{GrafakosH}  and  A. Seeger \cite {Seeger2} showed that $T_{a,1}$ is  bounded from $L^1$ to weak-$L^1$.

  S. Hofmann \cite{Hofmann2} proved  that $T_{\Omega ,a,k}$ is bounded from $L^p_w$ to itself for all $1<p<\infty$ and $w\in A_p$.

\par $\bullet$   {\bf Strongly singular multiplier operator}
$$T_{CF} f(x)= {\rm p.v.}\int\frac{e^{i|x-y|^{-b}}\chi_E(|x-y|)}{|x-y|^{n}} f(y)dy$$
where $0<b<\infty$ and $\chi_E$ is the characteristic function of the unit
interval $E=(0,1)\subset {\bf R} $.

$T_{CF}$ is bounded from  $L^p_w$ to itself for all $1<p<\infty$ and $w\in A_p$  and from $L^1_w$  to weak-$L^1_w$  for   $w\in A_1$, (see Chanillo \cite{Chanillo}  and see Fefferman \cite{F2} for the case $w=1$),
and from $H^1 _{w}$ into $L^1 _{w}$ for $  w\in A_1$,  (see \cite{Chanillo}  and see Fefferman and Stein \cite{FS72}  for the case $w=1$),
and  from  $H^1$ to itself   (see Sjolin \cite{Sjolin79}, see also Coifman \cite{Co1} for $n=1$).
    But $T_{CF}$ fails  to be bounded from   $H^p$ to itself for $0<p<1$ (see   \cite{FS72} or Sjolin \cite{Sjolin76}). And $T_{CF}$ fails to be bounded from $L^p _{|x|^{\alpha}}$  to itself for $1<p<\infty, \alpha\leq -n $ or $\alpha \geq n(p-1)$  (see \cite{Chanillo}).

 \par $\bullet$ {\bf  Partial sum operators of Fourier series }
 $$ C_{\xi}f(x)=\int_{{\bf R} }\frac{e^{2\pi i \xi y}}{x-y}    f(y)dy ,$$
 where $\xi\in {\bf R}$.

 The above bounded properties   of $H$ hold for $C_{\xi}$.

\par  $\bullet$ {\bf  Bochner-Riesz means at the  critical index}
 $$B_R^{(n-1)/2} f(x)=(f\ast K_R^{(n-1)/2})(x) $$
 with $K^{(n-1)/2}_R(x) =[(1-|\xi/R|^2)^{(n-1)/2}_+]~\breve{}~(x),$ where $\check{g}$ denotes the
inverse Fourier transform of $g$.

 $B_R^{(n-1)/2}$ is bounded from  $L^p_w$ to itself for all $1<p<\infty$ and $w\in A_p$  and from $L^1_w$  to weak-$L^1_w$  for   $w\in A_1$, ( see Shi and Sun \cite{SS}, and Vargas \cite{Vargas} and Christ \cite{Christ1} for $w=1$. But, $B_R^{(n-1)/2}$ fails to be  bounded  from  $L^1$  to itself.

 \par $\bullet$ {\bf  Oscillatory singular integrals operators}
$$T_{os}f(x)={\rm p.v.} \int_{{\bf
R}^n}e^{\lambda\Phi(x,y)}K(x,y)\varphi(x,y)f(y)dy,$$ where $K(x,y)$
is a Calder\'{o}n-Zygmund kernel,  $\varphi \in C_0^\infty({\bf
R}^n\times {\bf R}^n),$  $\lambda \in {\bf R}$, and $ \Phi(x,y)$ is
real-valued and  satisfies one of  the following conditions:

\par ~~ (a)  $\Phi(x,y)=(Bx,y)$ is a  real bilinear
form, $ \varphi=1$ and $\lambda=1$,

\par ~~ (b)  $\Phi(x,y)=P(x,y) $ is a polynomial, $ \varphi=1$ and $\lambda=1$,

\par ~~ (c) $\Phi(x,y)$ is a real analytic function on suup($\varphi $).

\par  $T_{os}$ is bounded from $L^p$ to itself for all $1<p<\infty$, (see, Phong and Stein \cite{PhongStein} for case (a), Ricci and Stein \cite{ RicciStein}) for case (b), and Pan \cite{Pan} for case (c)).

In the case (a), $T_{os}$ is bounded from $H_E^1$ to $L^1$, (see \cite{PhongStein}), but not from $H_E^p$ to $L^p$  for $0<p<1$, (see \cite{Pan3}),  where $H_E^1$ is an variant of the standard Hardy space $H^1$ and depends on $E=\Phi$.

   In the case (b), $T_{os}$ is bounded from $H_E^1$ to $L^1$, (see \cite{Pan3}), where $H_E^1$  depends on $E=\Phi$, but not from $H^1$ to $L^1$, (see \cite{PanY8}). It is bounded on $H_E^1$ if $\widetilde{T}(e^{\Phi(x,y)})=0$ for each $x\in {\bf R}^n$ in the BMO$_E$ sense, where $\widetilde{T}$ is the conjugate operator of $T$, (see
\cite{AlH} for the details). It is also bounded from $H^1$ to itself  if
$\Phi(x,y)=\Phi(x-y), K(x,y)=K(x-y)$ and $  \nabla \Phi (0)=0 $,
  (see \cite{HuP}).

  In the case (c),  $T_{os}$  is bounded from $H_E^1$ to $L^1$, (see \cite{PanY8}), where $H_E^1$  depends on $E=(\Phi,\lambda)$. It is bounded from $H^1$ to itself if
$\Phi(x,y)=\Phi(x-y), K(x,y)=K(x-y),\varphi(x,y)=\varphi(x-y), \Phi \in C_0^\infty, \nabla \Phi (0)=0 $ and $\partial^{\alpha}\Phi/\partial x^{\alpha}(0)\neq 0 $ for some multi-index $|\alpha|>1$,
  (see \cite{Pan4}).

The study of oscillatory singular integrals operators is also included in many other literatures, (see, for example,  \cite{FanP}).
There exists some oscillatory singular integrals operators  which  is bounded from
$L^p$ to itself for all $1<p<\infty$ but fails to be bounded from  $H^1$ to itself, (see \cite{Pan6}).

\par $\bullet$ {\bf  Oscillatory singular integrals operators  with rough kernel  }
$$T_{os,\Omega}f(x)={\rm p.v.} \int_{{\bf R}^n}e^{\lambda P(x,y)}\frac{\Omega (x-y) }{|x-y|^{n}}f(y)dy,$$
where  $ P(x,y) $ is a real-valued  polynomial, $\Omega$ satisfies (2.4) and (2.5).

By Lu and Zhang \cite{LuZ}, $T_{os,\Omega}$  is bounded from $L^p$ to itself for all $1<p<\infty$.
 By Challino and Christ \cite{ChanilloC},  this operator is bounded from $L^1$ to itself if $\Omega\in Lip(S^{n-1})$.

\par $\bullet$ {\bf  Pseudo-differential operator whose symbols in $ {\mathcal{S}}^{n(\varrho-1)}_{\varrho,\delta}$ with $0<\varrho\leq 1, 0\leq \delta<1$.}

 Let $m\in {\bf R},0<\varrho\leq 1, 0\leq \delta<1$.
A symbol in $ {\mathcal{S}}^{m}_{\varrho,\delta}$ will be a smooth function $p(x, \xi)$ defined on ${\bf R}^n\times {\bf R}^n$,  satisfying the estimates
$$
D^{\alpha}_xD^{\beta}_{\xi}p(x, \xi)\leq C_{\alpha,\beta}(1+|\xi|)^{m-\varrho |\beta| +\delta|\alpha|}.
$$
As usual, $\mathcal{L}^{m}_{\varrho,\delta}$
 denote the class of operators with symbol in $ {\mathcal{S}}^{m}_{\varrho,\delta}$.
When $0<\varrho\leq 1, 0\leq \delta<1$, and $m\leq (n+1)(\varrho-1)$, operators in $\mathcal{L}^{m}_{\varrho,\delta}$ are Calder\'{o}n-Zygmund operators, see Alvarez and  Hounie \cite {AH}.

\par For $m=n(\varrho-1)$, let $L\in \mathcal{L}^{n(\varrho-1)}_{\varrho,\delta}$. $L$ may not be in $\mathcal{L}^{m}_{\varrho,\delta}$ with  $0<\varrho\leq 1, 0\leq \delta<1$, and $m\leq (n+1)(\varrho-1)$, i.e. $L$ may not be
a Calder\'{o}n-Zygmund operators.
 Following \cite {AH},  $L$ has a distribution kernel $k(x, y)$  defined by
the oscillatory integral
 $$
 k(x,y)=(2\pi)^{-n} \int e^{i(x-y)\cdot \xi} p(x, \xi)d\xi,
 $$
which satisfies that
\begin{equation*}\label{2.1}
D^{\alpha}_x D^{\beta}_{y}k(x,y )\leq C|x-y|^{-n },~x\neq y.
\end{equation*}

By Theorem 3.4 and Theorem 3.2 in \cite {AH}, $L \in  {\mathcal{L}}^{n(\varrho-1)}_{\varrho,\delta}$ is bounded from $L^p$  to itself with $1<p<\infty$ and from $L^1$  to weak $L^1$,  since
$1\leq p<\infty, 0<\varrho\leq 1$ and $ 0\leq \delta<1$ imply $n(\varrho-1)\leq n(\varrho-1)|\frac{1}{p}-\frac{1}{2}|+ \min \{0,\frac{n(\rho-\delta)}{2}\}$.

$\varrho, \delta$ with $ 0<\varrho\leq 1$ and $ 0\leq \delta<1$  and $m=n(\varrho-1)$ satisfy the conditions of Theorem 5.2  and Theorem 5.4  in \cite {AH}, then, there is $p_0>0$ such that, for $p_0<p\leq 1$, $L \in  {\mathcal{L}}^{n(\varrho-1)}_{\varrho,\delta}$ is bounded from $H^p$ to $L^p$; moreover, if add $L^*(1)=0$ in the sense of $BMO$, then, $L  $ is bounded from $H^p$ to itself.

\par  By $T_1$ denote    the above linear  operators.
We have from Theorem 1.10 and Theorem 1.12 that

\begin{theorem}\label{th 8.1}
       Let $0<p <\infty, w\in A _{q,r}     $ with $ 0<q <p  $  and  $ 1<r <\infty  $.

\par     (i) If     $   1 < s< \infty  $ and $  rp/(r-1) \leq s   $, then each $T_1$  extends to a bounded operator from $BH^{p,s}_{w}$ to $L^{p}_{w}$.

 \par     (ii)   If     $\max\{ rp/(r-1), p/q\}< s< \infty$, then  each $T_1$ extends to a bounded operator from $BH^{p,s}_{w}$ to itself.

\end{theorem}

\par For each $T_1$, Theorem 2.2  extends  the $L^{p}_{}$ estimates with $1<p<\infty$   to all $0<p<\infty$. At the same time,  Theoorem 2.2 also gives a new weighted estimate.

\par  Theorem 2.2 is  sharp in the sense that the result may not hold for $ q= p =1$, since each $T_1$ fails to be bounded from $BH^{1,s}_{1}$ to $L^1$, noticing $BH^{1,s}_{1}=L^1$ for $s>1$ (see Proposition 3.7 below) and  $w=1\in  A _{1,r}   $ for  $ 1<r <\infty  $.

\par In particular, for $0<p \leq 1$, we have from Theorem 1.8 that

\begin{corollary}\label{th 8.1}
       Let $0<p \leq 1, w\in A _{q,r}     $ with $ 0<q <p  $  and  $ 1<r <\infty  $.  Then,
    \par  (i)   each $T_1$    extends to a bounded operator from $H^{p}_{w}$ to $L^{p}_{w}$,
    \par     (ii)   each $T_1$    extends to a bounded operator   from $H^{p}_{w}$ to itself.
\end{corollary}

\par For each $T_1$,  Corollary 2.3 gives the endpoint versions of the  $L^{p}_{}$ boundedness theorem  with $1<p<\infty$.

\subsection{Applications to the maximal operators }
Theorem 1.10(ii) and Theorem 1.12(ii) are applied to the following  maximal  operators.

\par $\bullet$ {\bf Maximal Hilbert transform}
   $$H^*f(x)=\sup_{\varepsilon> 0} |H^{\varepsilon}  f(x)|, $$
   where $H^{\varepsilon}f(x)= \int_{|x-y|>\varepsilon}\frac{f(y)}{x-y}dy  $, the truncated operator of $H$.

        $H^*$ is bounded from $L^p_w$ to itself for all $1<p<\infty$ and $w\in A_p$ and from $L^1_w$ to weak-$L^1_w$ for   $w\in A_1$,  (see Coifman and Fefferman \cite{CoFe}, and see Stein \cite{St61970} for the case $w=1$).

\par $\bullet$ {\bf Maximal  Riesz transform}
 $$R_j^*f (x)=\sup_{\varepsilon> 0} |R_j^{\varepsilon} f(x)|, j=1,2,\cdots,n,$$
where $R_j^{\varepsilon} f(x)=  \int_{|x-y|>\varepsilon}\frac{|x_j-y_j|}{|x-y|^{n+1}} f(y)dy $, the truncated operator of $R_j$.

 $R_j^*$ is bounded  from $L^p_w$ to itself for all $1<p<\infty$ and $w\in A_p$ and from $L^1_w$ to weak-$L^1_w$ for   $w\in A_1$,   (see Coifman and Fefferman \cite{CoFe}, and see Stein \cite{St61970} for the case $w=1$).

\par $\bullet$  {\bf Maximal Calder\'{o}n-Zygmund  operator}
$$T_{CZ}^*f(x)=\sup_{\varepsilon>0}|T_{CZ}^{\varepsilon}f(x)|$$
with
$$T_{CZ}^{\varepsilon}f(x)= \int_{|x-y|>\varepsilon} k(x,y)f(y)dy $$
for   $f\in L^p $ with compact support with $1<p<\infty$,  where $k(x,y)$   is a    Calder\'{o}n-Zygmund kernel   satisfying (2.2) and (2.3).

 $T_{CZ}^*$ is bounded  from $L^p_w$ to itself for all $1<p<\infty$ and $w\in A_p$ and from $L^1_w$ to weak-$L^1_w$ for   $w\in A_1$,   (see Meyer \cite{Mey}).

\par $\bullet$  {\bf Maximal   singular integral operator with rough kernel}
$$T^*_{\Omega} f(x)=\sup_{\varepsilon> 0} |T_{\Omega}^{\varepsilon} f(x)|,$$
with
$$T_{\Omega}^{\varepsilon} f(x)=  \int_{|x-y|>\varepsilon}\frac{\Omega (x-y)h(|x-y|)}{|x-y|^{n}} f(y)dy ,$$
 where $ h\in L^\infty ([0,\infty))$, and  $ \Omega $  satisfies (2,4) and (2.5)

 $T^*_{\Omega}$ is bounded from $L^p_w$ to itself for all $1<p<\infty$ and $w\in A_p$,  (see Duoandikoetxea and Rubio de Francia  \cite{DR}).

\par $\bullet$ {\bf Maximal strongly singular multiplier operator}
    $$T_{CF}^*f (x)=\sup_{\varepsilon> 0} |T_{CF}^{\varepsilon}\ast f(x)|.$$
 with
 $$T_{CF}^{\varepsilon} f(x)=  \int_{|x-y|>\varepsilon}\frac{e^{i|x-y|^{-b}}\chi_E(|x-y|)}{|x-y|^{n}} f(y)dy,$$
where $0<b<\infty$ and $\chi_E$ is the characteristic function of the unit
interval $E=(0,1)\subset {\bf R} $.

  $T_{CF}^*$ is bounded from $L^p_w$ to itself for all $1<p<\infty$ and $w\in A_p$,  (see Chanillo \cite{Chanillo}).

\par $\bullet$ {\bf Carleson operator}
$$ C^*f(x)=\sup_{\varepsilon>0}\sup_{\xi\in {\bf R}} \left|C_{\varepsilon, \xi}f(x)\right|$$
with
$$ C_{\varepsilon, \xi}f(x)=\int_{|x-y|>\varepsilon}\frac{e^{2\pi i \xi y}}{x-y}   f(y)dy,$$

$C^*$ is bounded from $L^p_w$ to itself for all $1<p<\infty$ and $w\in A_p$,  (see Hunt and Young \cite{HY}, and see Carleson \cite{Carl} and Hunt \cite{Hunt} for the case $w=1$).  But $C^*$ is not  bounded from $L^1$ to weak-$L^1$,  (see Kolmogorov \cite{Ko1}), and not even from $H^1$ to weak-$L^1$, (see Zygmund  \cite{Zy},  Chapter 8).

\par $\bullet$ {\bf Maximal  Bochner-Riesz means at the  critical index }
   $$B^*_{(n-1)/2} f(x)=\sup_{0<R<\infty}|B^R_{(n-1)/2} f(x) |.$$

  $B^*_{(n-1)/2}$ is bounded from $L^p_w$ to itself for all $1<p<\infty$ and $w\in A_p$  (see Shi and Sun \cite{SS}). But $B^*_{(n-1)/2}$ is not  bounded from $L^1$ to weak-$L^1$, and not even from $H^1$ to weak-$L^1$, (see Stein \cite{Stei}).

  \par $\bullet$ {\bf  Maximal   oscillatory singular integrals}
 $$ T_{P}^*f(x)=\sup_{\varepsilon>0}  \left|T_{P}^{\varepsilon}f(x)\right|$$
with
$$ T_{P}^{\varepsilon}f(x)=  \int_{|y|> \varepsilon} e^{  i P(x,y)} k(y)  f(x-y)dy,$$
where $P:{\bf R}^n\times {\bf R}^n\rightarrow {\bf R}$    is a polynomial
of two variables, and  $k$ is a suitable  Calder\'{o}n-Zygmund kernel on ${\bf R}^n$.
Krause and Lacey  proved   in  \cite{KrL } that $ T_{P}^*$ is bounded from $L^p_w$ to itself for all $1<p<\infty$ and $w\in A_1$.

\par $\bullet$ {\bf Polynomial Carleson operator}
$$ C_{d,n}f(x)=\sup_{P\in \mathfrak{P}_{d,n}}  \left| C_{P,d,n}f(x)\right|.$$
with
$$ C_{P,d,n}f(x)=\int_{{\bf R}^n} e^{  i P(x- y)} k(x-y)  f(y)dy,$$
where $\mathfrak{P}_{d,n}$ is the class of all real-coefficient polynomials in
$n$ variables with no constant term and of degree less than or equal to $d, d\in {\bf N}$, and  $k$ is a suitable  Calder\'{o}n-Zygmund kernel on ${\bf R}^n$.
Stein conjectured in \cite{Stein4} and  \cite{Stein5} that $ C_{d,n}$ is bounded on $L^p$ for any $1<p<\infty$, Lie proved the one dimensional case of this conjecture in  \cite{Lie}.

\par By $T_2$ denote    the  operators        $ H^{\varepsilon},R^{\varepsilon}_{j}, T^{\varepsilon}_{\Omega},T^{\varepsilon}_{F}, T^{\varepsilon}_{CZ}, T_{P}^{\varepsilon},  C_{\varepsilon, \xi} , B_{(n-1)/2 }$ and  $ C_{P,d,1}$ above,  by $T_2^*$ denote    the above maximal operators    associated with $T_2$.

It is known that each $T_2^*$  fails to be  bounded from $L^1$ to itself. At the same time, each $T_2^*$  fails to be  bounded from $H^p$ to itself for $0<p\leq 1$ since $H^p$ functions have the vanishing properties \cite{S3}.

  We have from Theorem 1.10 and Theorem 1.12 that

\begin{theorem}\label{th 8.1}
       Let $0<p <\infty, w\in A _{q,r}     $ with $ 0<q <p  $  and  $ 1<r <\infty  $.

\par     (i) If      $   1 < s<\infty  $ and $  rp/(r-1) \leq s   $, then each $T_2^*$  extends to a bounded operator from $BH^{p,s}_{w}$
to $L^{p}_{w}$.

 \par     (ii)   If $w\in    P  $   and  $\max\{ rp/(r-1), p/q\}< s< \infty$, then
each $T_2^*$  extends to a bounded operator from $BH^{p,s}_{w}$ to itself.

\end{theorem}

\par For each $T_2^*$, Theorem 2.4 extends the $L^{p}_{}$ estimates  with $1<p<\infty$  to all $0<p<\infty$. At the same time,  Theoorem 2.4 also gives a new weighted estimate.

\par  Theorem 2.4  is  sharp in the sense that the result may not hold for $ q= p =1$, since each $ T_2^*$  fails to be bounded from  $BH^{1,s}_{1}$ to $L^1$,  noticing $BH^{1,s}_{1}=L^1$ for $s>1$ (see Proposition 3.7 below) and  $w=1\in  A _{1,r}   $ for  $ 1<r <\infty  $.

\par
In particular, for $0<p \leq 1$,  we have from Theorem 1.8 that
\begin{corollary}\label{th 8.1}
       Let $0<p \leq 1, w\in A _{q,r}     $ with $ 0<q <p  $  and  $ 1<r <\infty  $.

\par     (i)      Each $T_2^*$   extends to a bounded operator from $H^{p}_{w}$ to $L^{p}_{w}$.

 \par     (ii)   If $w\in    P  $, then each $T_2^*$  extends to a bounded operator from $H^{p}_{w}$ to itself.

\end{corollary}

\par
 For each $T_2^*$,   Corollary 2.5  gives  the endpoint versions of the  $L^{p}_{}$ boundedness theorem  with $1<p<\infty$, including  Carleson-Hunt theorem.

For the polynomial Carleson operator $C_{d,n}$ of high dimensional,  we have

\begin{corollary}\label{th 8.1}
If  $C_{d,n}$ is bounded on $L^s$, then, Theorem 2.4 and Corollary 2.5 hold for $C_{d,n}$.

\end{corollary}

Corollary 2.6 gives a endpoint  version and a weighted version of the conjecture of Stein.

\par  By routine argument, the result for  Carleson operator in Theorem 2.4(i) implies:

\begin{corollary}\label{th 8.1}
       Let $0<p <\infty, w\in A _{q,r}     $ with $ 0<q <p  $  and  $ 1<r <\infty  $,   $   1 < s< \infty  $ and $  rp/(r-1) \leq s   $, then for all  $f\in BH^{p,s}_{w}$, we have
        \begin{equation}\label{2.2}
       \lim_{N\rightarrow \infty}\int_{|\xi|<N} \hat{f}(\xi)e^{2\pi ix\xi} d\xi=f(x),~~~w-a.e..
        \end{equation}
 In particular, for $0<p \leq 1, w\in A _{q,r}     $ with $ 0<q <p  $  and  $ 1<r <\infty  $, (2.7) holds for  all  $f\in H^{p}_{w}$.
\end{corollary}

Carleson-Hunt's theorem states that (2.7) holds for all  $f\in L^{p}_{}$ with $1<p<\infty$ and $w=1$, (see \cite{Carl,Hunt}). But it may not hold for $f\in L^{1}_{}$ (even $ H^{1}_{}$)(see \cite{Ko1,Zy}).
 Corollary 2.7 gives an extension to $0<p\leq 1$ for Carleson-Hunt's theorem.

\par Corollary 2.7 is sharp in the sense  that
the result may not hold for $ q= p =1$, since
there a function in  $BH^{1,s}_{1}$ such that (2.7) does not hold, (see \cite{Ko1} ), noticing
 $BH^{1,s}_{1}=L^1$ (see (3.9) below) and $1\in A _{1,r}$ with $1<r<\infty$.

The result for the maximal  Bochner-Riesz means at the  critical index in Theorem 2.4(i) implies:

\begin{corollary}\label{th 8.1}
     Let $n\geq 2$.   Let $0<p <\infty, w\in A _{q,r}     $ with $ 0<q <p  $  and  $ 1<r <\infty  $,     $   1 < s< \infty  $ and $  rp/(r-1) \leq s   $, then for all  $f\in BH^{p,s}_{w}$, we have
        \begin{equation}\label{ }
       \lim_{R\rightarrow \infty}\int (1-|\xi/R|^2)^{(n-1)/2}_+ \hat{f}(\xi)e^{2\pi ix\xi} d\xi=f(x),~~~w-a.e..
       \end{equation}\label{2.2}
In particular, for $0<p \leq 1, w\in A _{q,r}     $ with $ 0<q <p  $  and  $ 1<r <\infty  $, (2.8) holds for  all  $f\in H^{p}_{w}$.
\end{corollary}

It is known  that (2.8) holds for all  $f\in L^{p}_{}$ with $1<p<\infty$ and $w=1$,  (see  \cite{SS}). But it may not hold for $f\in L^{1}_{}$ (even $ H^{1}_{}$),  (see \cite{Stei}).
 Corollary 2.8 gives an extension to $0<p\leq 1$ for this result.

\par Corollary 2.8 is sharp in the sense that
the result may not hold for $ q= p =1$, since
there a function in  $BH^{1,s}_{1}$ such that (2.8) does not hold, (see \cite{Stei} ), noticing
 $BH^{1,s}_{1}=L^1$ (see (3.9) below) and $1\in A _{1,r}$ with $1<r<\infty$.

\par $\bullet$  {\bf   Hardy-Littlewood maximal operator  $M$}.

  $M$ is bounded from $L^p_w$ to itself for all $1<p\leq\infty$ and $w\in A_p$, (see Muckenhoupt  \cite{Muck}).

We have from Theorem 1.10 and Theorem 1.12 that

\begin{theorem}\label{th 8.1}
       Let $0<p <\infty. $
         \par (a) Let
        $ w\in D _{q}     \cap   P $ with $ 0<q <p  $,  then   $M$  extends to a bounded operator from $BH^{p,\infty}_{w}$
to $L^{p}_{w}$.

       \par (b) Let
        $ w\in A _{q,r}     \cap   P $ with $ 0<q <p  $  and  $ 1<r <\infty  $,

\par     (i) let     $   1 < s\leq \infty  $ and $  rp/(r-1) \leq s   $, then   $M$  extends to a bounded operator from $BH^{p,s}_{w}$
to $L^{p}_{w}$;

 \par     (ii)   let  $\max\{ rp/(r-1), p/q\}< s\leq  \infty$, then
  $M$  extends to a bounded operator from $BH^{p,s}_{w}$ to itself.

\end{theorem}

\par For $M$, Theorem 2.9 extends the $L^{p}_{}$    estimates with $1<p<\infty$   to all $0<p<\infty$.
 At the same time,  Theoorem 2.9 also gives a new weighted estimate.

\par In particular, for $0<p \leq 1$, we have from Theorem 1.8 that

\begin{corollary}\label{th 8.1}
       Let $0<p \leq 1$, and
            $ w\in A _{q,r}     \cap   P $ with $ 0<q <p  $  and  $ 1<r <\infty  $, then  $M$      extends to a bounded operator  from $H^{p}_{w}$ to $L^{p}_{w}$ and  from  $H^{p}_{w}$ to itself.
\end{corollary}

For  $M$,   Corollary 2.10  gives  the endpoint versions of the  $L^{p}_{}$ boundedness theorem  with $1<p<\infty$.

\par  Theorem 2.9  is  sharp in the sense that the result may not hold for $ q= p =1$, since $M$ fail to be bounded from $BH^{1,s}_{1}$ to $L^{1} $;
 Corollary  2.10 is  sharp in the sense that the result may not hold for $ q= p =1$, since $M$ fail to be bounded
 from $H^{1} $ to itself and  from $H^{1} $ to  $L^{1} $, noticing the facts $BH^{1,s}_{1}=L^1$ for $s>1$ (see (3.9) below),    $w=1\in A _{1,r}    \cap   P $ with $ 1<r <\infty  $, and $\{f\in L^1:Mf\in L^1\}=\{0\}$.

\par  Theorem 2.9(i) implies a Lebesgue differentiation theorem:
\begin{corollary}\label{th 8.1}
       Let $0<p <\infty, w\in A _{q,r}     \cap   P $ with $ 0<q <p  $  and  $ 1<r <\infty  $,      $   1 < s\leq \infty  $ and $  rp/(r-1) \leq s   $, then
      \begin{equation}\label{ }
       \lim_{x\in Q, l(Q)\rightarrow 0} \frac{1}{|Q|}\int_{Q} f(y) dy=f(x),~~~ a.e.,
       \end{equation}\label{ }
for all  $f\in BH^{p,s}_{w}$. In particular, for $0<p \leq 1, w\in A _{q,r}   \cap   P     $ with $ 0<q <p  $  and  $ 1<r <\infty  $, (2.9) holds for  all  $f\in H^{p}_{w}$.
\end{corollary}

 Corollary 2.9 gives an extension to $0<p< 1$ for the classical Lebesgue differentiation theorem that states that (2.9) holds for  all  $f\in L^{p}_{}$ with $1\leq p\leq  \infty$ (see \cite{St61970}).

\par Throughout the whole paper,  $C$ denotes a positive absolute constant not necessarily the same at each occurrence, and a subscript is added when we wish to make clear its dependence on the parameter in the subscript.
  For any cube $Q$ and $\lambda>0$,  $ \lambda Q $ denotes the cube concentric with $Q$ whose each edge is $\lambda$ times as long.
For $1\leq s\leq \infty$, $s'$ denotes the conjugate of $s$, which satisfies $1/s+1/s'=1$. By $A_p$  with  $1\leq p\leq \infty$ we denot the classical  Muckenhopt  class.

We express our gratitude   to   David Cruz-Uribe   for his comments  on our other article
which led to substantial  improvements of this paper.

\section
 {Some elementary results }
\par In this section we collect a few standard facts and elementary results for weights  and $B^{p,s}_w$ that
will be used later.

\begin{lemma}\label{pro_2.9}
  Let $1<r< \infty$.  We have ,
  \par (i) if  $w\in RH_{r}$, then a.e. implies $w$-a.e.,

  \par (ii) if  $w\in RH_{r}\bigcap  P $, then $w$-a.e. implies a.e..

\end{lemma}

\par \proof

 It is known that $w\in A_{\infty}$  if and only if $w \in  RH_r$ for some $r > 1$, and that if $w\in A_{\infty}$ then there is a $1<p_0<\infty$ such that  $w\in A_{p_0}$. It follows, for $w \in  RH_r$, there is a $1<p_0<\infty$ such that
\begin{equation}\label{2.3}
 C_1
\left( \frac{|E|}{|Q|}
\right )^{p_0}
\leq
 \frac{w(E)}{w(Q)}\leq C_2
\left( \frac{|E|}{|Q|}
\right )^{(r-1)/r}
\end{equation}
for any measurable subset $E$ of a cube $Q$, where $C_1, C_2>0$ are    constants independing on $Q$  and $E$, (see \cite{GRubio}).

For (i), we need to prove that  $w(E)=0$ if $|E|=0$ for a measurable set $E$.
To do this, we set ${\bf R}^n=\bigcup_{i=1}^{\infty} Q_i$, where each $Q_i$ is a cube with $|Q_i|=1$.
By $w \in  RH_r$, (3.1) holds for each $Q_i$ and any measurable subset $E_i$ of  $Q_i$ with $|E_i|\neq 0$, the left inequality of (3.1) implies
  $w(Q_i)\neq \infty$ for each $i$. At the same time,   if $|E|=0$, clearly $|Q_i\bigcap E|=0$, then  we have from the right inequality of (3.1) that $w(Q_i\bigcap E)=0$ for each $i$,  it follows that  $w( E)\leq \sum w(Q_i\bigcap E)=0$.
\par
For (ii),  we need to prove that $|E|=0$ if $w(E)=0$.  By $w\in  P $, then, there exist a sequence $\{Q_i\}$ of cubes whose interiors are disjoint each other such that  ${\bf R}^n=\bigcup_{i=1}^{\infty} Q_i$ and $w(Q_i)>0$ for each $i$. By  $w\in RH_{r}$,  then, (3.1) holds
for each $Q_i$ and  any measurable subset $E_i$ of  $Q_i$ with $w(E_i)\neq 0$, the right inequality of (3.1) implies
 $|Q_i|\neq \infty$ for each $i$.
Let $w( E)=0$,  by the left inequality of (3.1), we have
\begin{equation*}\label{2.1}
\left(
\frac{|Q_i\bigcap E|}{|Q_i|}
\right )^{p_0}
\leq C \frac{w(Q_i\bigcap E)}{w(Q_i)} = C \frac{0}{w(Q_i)} =0
 \end{equation*}
for each $i$, it follows $|Q_i\bigcap E|=0$, then, $| E|\leq \sum |Q_i\bigcap E|=0$.

 \par The lemma have been proved.

The results in Lemma 3.1 hold when $ RH_{r}$ is replaced by $A_{\infty}$.

\begin{lemma}\label{lem_3}
Let $0<p <\infty ,  w\in  RH_r  $ with   $ 1<r <\infty  $,  and $  rp/(r-1) \leq s \leq \infty $. Then, for any cube $ Q$ in ${\bf R}^n$,
\begin{equation}\label{2.2}
\left(\int_{Q}w^{ s/(s-p)}\right)^{(s-p)/s} \leq C |Q|^{-p/s} w(Q ).
 \end{equation}
\end{lemma}

\par \proof
For $s=\infty$, noticing that $(s-p)/s=1$, (3.2) is obvious. For $s<\infty$,
by H\"{o}lder inequality and (1.2), we have
\begin{eqnarray*}
\left(\int_{Q}w^{ s/(s-p)}\right)^{(s-p)/s}
&\leq &
 |Q|^{(s-p)/s}  \left(\frac{1}{|Q|}\int_{Q}w^{ r}\right)^{1/r}
      \\
   &\leq & C |Q|^{-p/s}   w(Q ).
 \end{eqnarray*}
Thus, Lemma 3.2 holds.

For  $0<p<\infty$ and weigh $w$, we have
\begin{equation}\label{2.1}
\|\sum_{i=1}^{\infty}  a_i \|^{\bar{p}}_{L^{p}_{w}} \leq  \sum_{i=1}^{\infty}\|  a_i \|^{\bar{p}}_{L^{p}_{w}}
\end{equation}
for $a_i\in L^p_w, i=1,2,3,\cdots$, that follows from Minkowski inequality for $1\leq p\leq \infty$ and the inequality $(|a|+|b|)^p\leq |a|^p+|b|^p$ for $0<p\leq 1$.

\begin{proposition}\label{Th3.1}
 Let $ 0< p <  \infty$.
  Let $ w\in RH_{r}$ with $1<r<\infty$ and  $rp/(r-1)\leq s \leq \infty $ or  $w$ be a weight and  $  s=\infty$.
 Let $\sum_{i=1}^{\infty} \lambda _i a_i  \in B^{p,s}_{w}$, where
  each $a_i $ is a $(p,s,w)$-block  and $\sum_{i=1}^{\infty}|\lambda_i|^{\bar{p}} <\infty$.
Then,
\par (i) for every $a_i$,
\begin{equation}\label{3.1}
 \|a_i\|_{L^{p}_{w}} \leq C;
  \end{equation}

 \par (ii)   $\sum_{i=1}^{\infty}\lambda_i a_i$ converge in $ L^{p}_{w} $ and $w$-a.e., and
 \begin{equation}\label{3.1}
 \|\sum_{i=1}^{\infty} \lambda _i a_i \|_{L^{p}_{w}} \leq C \left(\sum_{i=1}^{\infty}|\lambda_i|^{\bar{p}}\right)^{1/{\bar{p}}}.
  \end{equation}
    Consequently, $BL^{p,s}_{w} \subset L^{p}_{w}$.
\end{proposition}
\par \proof
 \par (i). Let $a_i$ be a $(p,s,w)$-block, and suppose that $\supp a_i \subset Q$ with the cental  $x_0$.  For $rp/(r-1)\leq s < \infty $,   we see from $r>1$  that $p<s$ and  $s/(s-p)<r$, then,
using the H\"{o}lder inequality, (3.2) and the definition of $a_i$,
 we have
\begin{eqnarray*}
\|a_i\|^p_{L^{p}_{w}}
  &\leq & \|a_i\|^p_{L^{s}}
    \left(\int_{Q
}w^{ s/(s-p)}dx\right)^{(s-p)/s}
  \\ &\leq & C \|a_i\|^p_{L^{s}}
 |Q|^{-p/s}
  w (Q)
\\
& \leq & C.
\end {eqnarray*}

\par For $s=\infty,$ we have
$
\|a_i\|^p_{L^{p}_{w}}
  \leq  \|a_i\|^p_{L^{\infty}}  \int_{Q
}wdx \leq  1
$.
Thus, (3.4) holds,  (i) have been proved.

 (ii). First, using (3.3), (3.4)  and  $\sum_{i=1}^{\infty}|\lambda_i|^{\bar{p}} <\infty$, we see that  $  \{ \sum_{i=1}^{N} \lambda _i a_i \}_{N=1}^{\infty} $ is a Cauchy sequence   in $L^{p}_{w}$, and  by  the completeness of $L^{p}_{w}$,  there is a unique $f\in L^{p}_{w}$, such that
 $   \sum_{i=1}^{N} \lambda _i a_i \rightarrow f $ in $L^{p}_{w}$ as $N\rightarrow\infty$.
\par  And then, we denote $   f=\sum_{i=1}^{\infty} \lambda _i a_i $, using (3.3) two times   and (3.4) again,
   we have that
 \begin{eqnarray*}\label{3.1}
 \|\sum_{i=1}^{\infty} \lambda _i a_i \|^{\bar{p}}_{L^{p}_{w}}
 &=&  \|f \|^{\bar{p}}_{L^{p}_{w}}
 \leq  \|\sum_{i=1}^{N} \lambda _i a_i \|^{\bar{p}}_{L^{p}_{w}} + \|f-\sum_{i=1}^{N} \lambda _i a_i \|^{\bar{p}}_{L^{p}_{w}}
 \\
 &\leq &  C  \sum_{i=1}^{N}|\lambda_i|^{\bar{p}} + \|f-\sum_{i=1}^{N} \lambda _i a_i \|^{\bar{p}}_{L^{p}_{w}}
  \end{eqnarray*}
   for all $N>1$,
letting $N\rightarrow \infty$, we get (3.5), at the same time, we have $B^{p,s}_{w} \subset L^{p}_{w}$, it follows $BL^{p,s}_{w} \subset L^{p}_{w}$ by the definition.

\par Finally,
 we need to prove that  $\sum_{i=1}^{N} \lambda _i a_i$ converges also to $ f $   $w$-a.e.. To do this, we need only to prove for all $\delta>0$ that
\begin{equation}
w(\{ x:|f-\sum_{i=1}^{\infty}\lambda _ia_i|>\delta\})=0.
\end{equation}
In fact, by using   (3.3), (3.4),    $\sum_{i=1}^{\infty}|\lambda_i|^{\bar{p}}<\infty$, and the fact that $\sum_{i=1}^{N} \lambda _i a_i$ converges  to $ f $   in $L^{p}_{w}$,
we have for a given $\delta>0$ and  all $N>1$ that
\begin{eqnarray*}
w(\{ x:|f-\sum_{i=1}^{\infty}\lambda _ia_i|>\delta\})
&\leq &
\int_{\{ x:|f-\sum_{i=1}^{N}\lambda _ia_i|>\delta/2\}} w
+\int_{\{ x:|\sum_{i=N+1}^{\infty}\lambda _ia_i|>\delta/2\}} w
\\
&\leq &
\left( \frac{2}{\delta}\right)^p
\left(\int|f-\sum_{i=1}^{N}\lambda _ia_i|^p w
+
\int|\sum_{i=N+1}^{\infty}\lambda _ia_i|^pw\right)
\\
&\leq &
\left( \frac{2}{\delta}\right)^p
\left(\int|f-\sum_{i=1}^{N}\lambda _ia_i|^p w
+
\left( \sum_{i=N+1}^{\infty}|\lambda _i|^{\bar{p}}\right)^{p/\bar{p}}\right)
,
\end{eqnarray*}
letting $N\rightarrow \infty$, we get (3.6),  which implies $\sum_{i=1}^{N} \lambda _i a_i$ converges  to $ f $   $w$-a.e..
  Then,  (ii) have been proved. Thus, we finish the proof of Proposition 3.3.

\begin{proposition}\label{pro_3.14}
Let $   p , s ,w$ as in Definition 1.4.
Let $ g$ is measurable, and  $|g(x)|\leq |f(x)|$, a.e..
    If   $f\in B^{p,s}_{w}$,  then $g\in B^{p,s}_{w}$ and
\begin{equation*}\label{3.9}
\|g\|_{B^{p,s}_{w}}\leq \|f\|_{B^{p,s}_{w}}.
\end{equation*}

\end{proposition}

\par \proof  The proof of the proposition is similar to that of Proposition 2.11 in \cite{LTW1}.

\begin{proposition}\label{pro_3.11}
  Let $   p , s ,w$ as in Definition 1.4.
 If $f_i(x) \in B^{p,s}_{w}, i=1,2,\cdots$, and $\sum_{i=1}^{\infty} \|f_i\|^{\bar{p}}_{B^{p,s}_{w}}<\infty$, then $\sum_{i=1}^{\infty} f_i \in {B^{p,s}_{w}}$ and
 \begin{equation}\label{3.8}
\|\sum_{i=1}^{\infty} f_i\|^{\bar{p} }_{B^{p,s}_{w}}\leq \sum_{i=1}^{\infty} \|f_i\|^{\bar{p}}_{B^{p,s}_{w}}.
 \end{equation}

\end{proposition}

\par \proof     Let $f_i \in B^{p,s}_{w}, i=1,2,\cdots.$ For
any $\varepsilon >0$ and each $f_i$, there exists
a sequence  $\{b^{(i)}_k\} $  of $(p,s,w)$-blocks  and a sequences $\{m^{(i)}_k\} $   of real numbers with $\sum_{k=1}^{\infty}|m^{(i)}_k|^{\bar{p}} <\infty$, such that
 $f_i=\sum_k m^{(i)}_kb^{(i)}_k $
   and
$\|f_i\|_{B^{p,s}_{w}}^{\bar{p}} \geq
 \sum_k|m^{(i)}_k|^{\bar{p}}  -\frac{1}{2^i}\varepsilon. $
 It follows
  $\sum_if_i=\sum_i\sum_k m^{(i)}_kb^{(i)}_k $
   and
 $
 \sum_i \sum_k|m^{(i)}_k|^{\bar{p}}\leq
\sum_i \|f_i\|_{B^{p,s}_{w}}^{\bar{p}}+\varepsilon \sum_i\frac{1}{2^i} <\infty.
 $
Thus, by the definition of $B^{p,s}_{w}$, we have  $\sum_if_i\in B^{p,s}_{w}$, and
$$\|\sum_if_i\|_{B^{p,s}_{w}}^{\bar{p}} \leq
 \sum_i \sum_k|m^{(i)}_k|^{\bar{p}}\leq
\sum_i \|f_i\|_{B^{p,s}_{w}}^{\bar{p}}+\varepsilon \sum_i\frac{1}{2^i}, $$
letting  $\varepsilon \rightarrow 0$,
(3.7) follows. Thus, we have proved the proposition.

\begin{proposition}\label{pro_3.16}

 Let $   p , s ,w$ as in Definition 1.4, then
   $ B^{p,s}_{w} $  is  complete.
\end{proposition}

\par \proof
Let $\{u_k\}$ be a Cauchy sequence in $ B^{p,s}_{w}. $
 For any $\varepsilon >0$, there exists a subsequence $\{u_{k_j}\} $
  of $\{u_k\}$ such that
  $$\|u_{k_{j+1}}-u_{k_j}\|_{B^{p,s}_{w}} < \frac{\varepsilon}{2^j}, ~~~~j=1,2,\cdots, $$
and
$$\|u_{k }-u_{k_1}\|_{B^{p,s}_{w}} <
\varepsilon $$ for $k>k_1.$
 Set
$\bar{u}=\sum_{j=1}^{\infty}(u_{k_{j+1}}-u_{k_j})$.  Since $u_{k_{j+1}}-u_{k_j}\in B^{p,s}_{w}$ for $j=1,2,\cdots$, and
$$
\sum_{j=1}^{\infty}\|(u_{k_{j+1}}-u_{k_j})\|_{B^{p,s}_{w}}
< \sum_{j=1}^{\infty} \frac{\varepsilon}{2^j}=\varepsilon<\infty,
$$
 we have by Proposition 3.5 that
 $\bar{u}\in B^{p,s}_{w}$ and
 $$ \|\bar{u}\|_{B^{p,s}_{w}}
<\varepsilon.$$
 Let
$u=\bar{u}+u_{k_1}$,
using Proposition 3.5 again, we have
 $u\in B^{p,s}_{w},$ and
 $$\|u_{k }-u_{ }\|_{B^{p,s}_{w}} \leq \|u_{k }-u_{k_1 }\|_{B^{p,s}_{w}} +
\|\bar{u} \|_{B^{p,s}_{w}}
 < 2
\varepsilon $$ for $k>k_1.$ Thus, we have proved the proposition.

\begin{proposition}\label{pro_3.13}
  Let $0< p < \infty, 0<s_1\leq s_2\leq \infty, w$ be a weight,
 then
\begin{equation}
B^{p,s_2}_{w} \subset B^{p,s_1}_{w}.
\end{equation}
In particular, for $1< s\leq \infty $ and $w=1$, we have
\begin{equation}
L^1=  B^{1,s}_{1}.
\end{equation}

\end{proposition}

\proof By the H\"{o}lder inequality, we see that a $(p,s_2,w)$-block is a $(p,s_1,w)$-block  for $0<s_1\leq s_2\leq \infty$, (3.8) follows.

We see ${ B}^{1,\infty}_1 =L^1 $  from \cite{S3}, pages 112 and  129, by (3.8), it follows
$ L^1 \subset B^{1,s}_{1}$ for $0< s\leq \infty $.
  Noticing    $w=1 \in RH_r$  for $1<r<\infty$,
  by  Proposition 3.3(ii), we have
$  B^{1,s}_{1}\subset L^1$ for $1 < s \leq \infty $.
(3.9) follows. Thus, the proposition holds.

\begin{lemma}\label{pro 7.4}
  Let $0<p <\infty. $  Suppose that an operator  $ T $   has definition for a $(p,s,w)$-bolck $h$ and  satisfies (1.5) and (1.6).
  If
  \par
  (a) $ w\in D _{q}     $ with $ 0<q <p  $,  and $ s = \infty $, or
  \par
  (b)
  $ w\in A _{q,r}     $ with $ 0<q <p  $  and  $ 1<r <\infty  $,  and $ \max \{ rp/(r-1),1\}\leq s \leq \infty$,
 \\ then
\begin{equation}\label{7.5}
\|T h\|_{L^{p}_{w}}\leq
C.
\end{equation}
\end{lemma}

\par \proof
Let  $h$ be a $(p,s,w)$-block. To prove (3.10), we suppose that supp $h\subseteq Q$,
    and  write
\begin{eqnarray*}
 \|Th\|^p_{L^p_w} =\int_{2n^{1/2}Q}|Th|^pw+\int_{{\bf R}^n \backslash 2n^{1/2}Q}|Th|^pw=: I+II.
 \end{eqnarray*}
For $I$,  when $s=\infty$, we have  by (1.6), (1.1) and the definition of $h$ that
$$I\leq C\|h\|^p_{L^{\infty}}w(2n^{1/2}Q)
\leq C, $$
 when $s<\infty$,
\begin{eqnarray*}
I
&\leq&
\left(\int_{2n^{1/2}Q}|Th|^s\right)^{p/s}
\left(\int_{2n^{1/2}Q}w^{s/(s-p)}\right)^{(s-p)/s}
\\
&&
~~~~~~~~~~~~~~~~({\rm by~H\ddot{o}lder~inequality~ for ~the~ index} ~s/p)
\\
&\leq&
C \left(\int_{Q}|h|^s\right)^{p/s}
 |Q|^{-p/s}w(Q)^{}
\\
&&~~~~~~~~~({\rm by}~(1.6),(3.2) {\rm ~and~} (1.1))
\\
&\leq &
C
\\
&&~~~~~~~~~~~~~~~~~~~~({\rm by~the~definition~ of ~}h).
 \end{eqnarray*}

 To estimate  $II$, we suppose that $Q$ has the center $x_0$,
we have
\begin{eqnarray*}
II
 &=&
\int_{{\bf R}^n\backslash 2n^{1/2}Q}|Th|^pw
\\
&\leq &
\|h\|_{L^1}^p
\sum_{i=1}^{\infty}
\int_{{2^{i+1}n^{1/2}Q\backslash 2^in^{1/2}Q}}|x-x_0|^{-np}w(x)dx
\\
&&~~~~~~~~~~~~~~~~~~~~~~~~~~~~~~~~~~~({\rm by}~(1.5))
\\
&\leq &
C \|h\|_{L^1}^p
\sum_{i=1}^{\infty} 2^{-npi}|Q|^{-p}w({2^{i+1}n^{1/2}Q})
\\
&\leq &
C \|h\|_{L^1}^p|Q|^{-p}w(Q)
\sum_{i=1}^{\infty} 2^{-n(p-q)i}
\\
&&~~~~~~~~~~~~~~~~~~~~~~~({\rm by~(1.1)~since}~~w\in D_{q})
\\
&\leq &
C \|h\|_{L^s}^p|Q|^{p/s'}|Q|^{-p}w(Q)
\\
&&~~~~~~~({\rm since} ~q<p~ {\rm and ~H\ddot{o}lder~inequality})
\\
&\leq &
C
\\
&&~~~~~~~~~~~~~~~~~~~~~~~~~~~~~~~~~~~~~~~~~~~~~({\rm by~the~definition~ of ~}h).
 \end{eqnarray*}
\par    The proof of Lemma 3.8 is complete.

\section
 { Block characterization of  weighted Hardy spaces }

\par In this subsection, we  prove Theorem 1.8. It follows from the following two theorems.

\begin{theorem}\label{Th3.1}
 Let $ 0< p <  \infty$. Let $\sum_{i=1}^{\infty} \lambda _i a_i  \in B^{p,s}_{w}$, where
  each $a_i $ is a $(p,s,w)$-block  and $\sum_{i=1}^{\infty}|\lambda_i|^{\bar{p}} <\infty$.
   Suppose that
  \par
  (a) $ w\in D _{q}     $ with $ 0<q <p  $,  and $ s = \infty $, or
  \par
  (b) $ w\in A _{q,r}      $ with $ 0<q <p  $  and  $ 1<r <\infty  $,
and  $ 1< s\leq \infty  $ and $ rp/(r-1)\leq s $.
Then,

\par (i) for every $a_i$,
 \begin{equation}\label{3.1}
 \|a_i\|_{H^{p}_{w}} \leq C,
  \end{equation}
\par (ii)
$\sum_{i=1}^{\infty}\lambda_i a_i$ converge in $ H^{p}_{w} $ and $w$-a.e., and
\begin{equation}\label{3.1}
 \|\sum_{i=1}^{\infty} \lambda _i a_i \|_{H^{p}_{w}} \leq C \left(\sum_{i=1}^{\infty}|\lambda_i|^{\bar{p}}\right)^{1/{\bar{p}}}.
  \end{equation}
  Consequently, $BH^{p,s}_{w} \subset H^{p}_{w}$.
\end{theorem}

 \proof
\par (i). To prove (4.1), by the definition of $H^p_w$,  we need only to prove that, for a $(p,s,w)$-block $h$, there a constant $C$ independent of $h$ such that
 \begin{equation}\label{ }
\|M_{\varphi}h\|_{L^p_w} \leq C
\end{equation}
for some $\varphi\in \mathcal{S}$ with  $\int_{{\bf R}^n}\varphi (x) dx=1$. To do this, we take $\varphi\in \mathcal{S}$ such that supp $\varphi \subset \{x:|x|\leq 1\}$ and $\int_{{\bf R}^n}\varphi (x) dx=1$,  and suppose that $\supp h \subset Q$ with the cental  $x_0$. Once it is true
 that $M_{\varphi}$ satisfies (1.5)  for $x\in (2n^{1/2}Q)^c $ and (1.6),  then,    we have  from Lemma 3.8 that (4.3) holds.

\par   In fact, when $x\in 2n^{1/2}Q, M_{\varphi}h(x)$ is controlled by Hardy-Littlewood maximal function $Mh(x)$, by the $L^s$ boundedness of $M$ with $1<s\leq \infty$, $Mh$ satisfis (1.6),
  $M_{\varphi}h$ follows.

 When $x\in (2n^{1/2}Q)^c $, for $y\in Q$,   we see that $|x-y|\geq
l(Q)/2\geq |y-x_0| ,$  it follows that $|x-x_0|\leq |x-y| + |y-x_0|\leq 2|x-y|$.
At the same time,  we have by $\varphi\in \mathcal{S}$ that $|x|^{n}|\varphi(x)| \leq C_n\sum _{\beta=n} |x^{\beta}||\varphi(x)|\leq C<\infty $, where $\beta$ are multiindex,
 see page 95 in \cite {Grafakos},
 it follows
 $|\varphi(x)| \leq C_n|x|^{-n} $.
Then, for $x\in (2n^{1/2}Q)^c $,
\begin{eqnarray*}
|\varphi_t\ast h(x)|=|\int_Q\frac{1}{t^n}\varphi(\frac{x-y}{t})h(y)dy|\leq C \frac{1}{|x-y|^n} \|h\|_{L^1}\leq C \frac{1}{|x-x_0|^n} \|h\|_{L^1},
\end{eqnarray*}
it follows that  $M_{\varphi}h$ satisfies (1.5) for $x\in (2n^{1/2}Q)^c $.  Thus, we have proved (i).

\par (ii). Using the known inequality $\|f+g\|^{\bar{p}}_{H^p_w}\leq \|f\|^{\bar{p}}_{H^p_w}+\|g\|^{\bar{p}}_{H^p_w}$ ,   (4.1) and  $\sum_{i=1}^{\infty}|\lambda_i|^{\bar{p}} <\infty$, we see that  $  \{ \sum_{i=1}^{N} \lambda _i a_i \}_{N=1}^{\infty} $ is a Cauchy sequence   in $H^{p}_{w}$, and then, by  the completeness of $H^{p}_{w}$,  there is a unique $f\in H^{p}_{w}$, such that
 $   \sum_{i=1}^{N} \lambda _i a_i \rightarrow f $ in $H^{p}_{w}$ as $N\rightarrow\infty$.
At the same time, it follows that,  for any $\varepsilon>0$, there is a $N>0$ such that
 $\|f  -\sum_{i=1}^{N}a_i ) \|_{H^{p}_{w}}<\varepsilon $, we then have
\begin{eqnarray*}
\|f\|^{\bar{p}}_{H^{p}_{w}} \leq \sum_{i=1}^{N}|\lambda_i|^{\bar{p}}\|a_i \|^{\bar{p}}_{H^{p}_{w}} +\|f  -\sum_{i=1}^{N}a_i ) \|^{\bar{p}}_{H^{p}_{w}}
 \leq C\sum_{i=1}^{N}|\lambda_i|^{\bar{p}} +\varepsilon^{\bar{p}}
\end{eqnarray*}
which implies (4.2). And from Proposition 3.3, we have that $\sum_{i=1}^{\infty}a_i $ converge $w$-a.e..  Thus, we have proved (ii).

The proof of Theorem 4.1 is complete.

Next, we want to prove $H^{p}_{w} \subset BH^{p,s}_{w}$, to do this, we  recall that we can also characterize  $H_w^p$ in terms of
atoms in the following way (see \cite{Gc}). Let $0<p\leq 1\leq s \leq \infty, p\neq s$ and $ w\in A_s$ with the critical index $\tilde{s}_w=\inf\{s>1:w\in A_s\}$, $N\geq N_0=[n(\tilde{s}_w/p-1)]$,
  a function $a$ is said a  $w$-$(p,s,N)$-atoms if
    \par (i)~~~~supp $ a\subseteq Q,$  a cube in ${\bf R}^n$,
\par (ii)~~~~$\|a\|_{L_w^{s}   }\leq w(Q)^{1/s-1/p},$
\par (iii)~~~~$\int_{{\bf R}^n} a(x)  x^{\alpha} dx =0 $ for every multi-index $\alpha$ with $|\alpha|\leq N$.

\par  Let $H^{p,s,N_0}_w$ denote the space consisting of all
$\sum \lambda_ia_i$
 that converge in $H^p_w$,
 where each $a_i$ is a $w$-$(p,s,N_0)$-atom and $\sum |\lambda_i|^p<\infty.$

\begin{theorem}\label{Th3.1}
  Let $ 0< p \leq 1, w\in  A_{\infty}$,
   and $1\leq s \leq \infty$,
then,
\begin{equation}\label{3.1}
H^{p}_{w} \subset BH^{p,s}_{w} .
\end{equation}
\end{theorem}

\par \proof  By the atom characterization of $H_w^p$  in \cite{Gc},
   for $0<p\leq 1$ and $w\in A_{\infty}$,
$$H^p_w=H^{p,\infty,N_0}_w. $$
For every $N\geq 0$, we see from the definitions that each $w$-$(p,\infty,N)$-atom is a $(p,\infty,w)$-block for $0<p\leq 1$ and $w\in A_{\infty}$, it follows
that
$$H^{p,\infty,N_0}_w\subset BH^{p,\infty}_w.$$
On the other hand,   we see by (3.8)that
$$BH^{p,\infty}_w\subset BH^{p,s}_w$$
for $1\leq s\leq \infty$.  (4.4) follows. The proof of Theorem 4.2 is complete.

\par
\proof [Proof of Theorem 1.8]
Noticing $RH_r\subset A_{\infty}$ for $1<r<\infty$,  (1.3) follows from  the above Theorem 4.1 and Theorem 4.2. On the other hand, it is well known  that $H^1\neq L^1$ which is equal to $  B^{1,s}_{1}$ by (3.9), i.e. (1.3) does  not hold when $ q= p =1$ and $w=1$ which is in $A _{1,r}   $ for $1<r<\infty$. Thus, we have proved Theorem 1.8.

\par
\section
 {
   Hardy-Littlewood maximal function characterization of   weighted Hardy spaces
  }

In this section, we  prove Theorem 1.9. We begin with the following fact.

\par Let $ 0< p <  \infty$. Let $w\in RH_{r}\cap P $ with $1<r<\infty$ and   $\max\{rp/(r-1),1\}\leq  s \leq \infty $, or,
 $w\in  P $ and  $  s=\infty$.
 Let $f(x)=\sum_{l}\lambda_l a_l\in  BL^{p,s}_w$, where $\{a_l\}$ is a sequence   $(p,s,w)$-blocks and $\{\lambda_l\}$ is a sequence real numbers with $ \sum_{l}|\lambda_l|^{\bar{p}} <\infty$.

By Proposition 3.3, we see that $f(x)=\sum_{l}\lambda_l a_l$ converges   $w$-a.e., by Lemma 3.1, it holds   a.e., then,
\begin{equation*}\label{ }
|f(x)|\leq \sum_{l}|\lambda_l| |a_l(x)|,~a.e..
 \end{equation*}
Each $Ma_l$ is well defined since $a_l\in L^s$ with $1\leq s\leq \infty$.  Then,
we have by   Minkowski inequality   that
 \begin{equation}\label{ }
Mf(x)\leq \sum_{l}|\lambda_l| Ma_l(x)
 \end{equation}
holds a.e..

\par Theorem 1.9 is contained in the following two theorems.

\begin{theorem}\label{th_4.5}
    Let $0<p <\infty  , w\in A _{q,r}     \cap   P $ with $ 0<q <p  $  and  $ 1<r <\infty  $,
   $  1 < s\leq \infty  $  and  $  rp/(r-1) \leq s   $.
   If
   $f(x)=\sum_{l}\lambda_l a_l\in BL^{p,s}_{w},
 $
     where each $\{a_l\}$ is a $(p,s,w)$-block and $ \sum_{l}|\lambda_l|^{\bar{p}} <\infty,$ then $f\in ML^p_{w}$
  and
   \begin{equation}\label{4.2}
 \|f\|_{ML^p_{w}} \leq C\left(\sum_{l}|\lambda_l|^{\bar{p}}\right)^{1/\bar{p}}.
 \end{equation}
 Consequently, $BH^{p,s}_{w} \subset ML^{p}_{w}$.
\end{theorem}

\par \proof   Let $f(x)=\sum_{l}\lambda_l a_l \in BL^{p,s}_{w} $.
    By Lemma 3.1, (5.1) holds $w$-a.e., then we have by using (3.3)  that
\begin{equation}\label{4.2}
 \|Mf\|^{\bar{p}}_{L^p_{w}} \leq \sum_{l}|\lambda_l|^{\bar{p}}  \|Ma_l\|^{\bar{p}}_{L^p_{w}}.
 \end{equation}
Once
 \begin{equation}\label{ }
\|Ma_l\|_{L^p_w} \leq C
\end{equation}
is proved for a constant $C$  independent  of $a_l$,  (5.2) follows from (5.3). To prove (5.4), by Lemma 3.8, we need only to check that
 $M$ satisfies (1.5) and (1.6).
 To do this, we suppose that $\supp a_l\subset Q_l$, a cube with centre $x_0$. Let $x\in (2n^{1/2}Q_l)^c$, for any $y\in Q_l$,    we have   $|x-y|\geq |x-x_0|/2$, then for a cube $Q$ containing $ x$ and with side length $l(Q)\leq |x-x_0|/2$, we have that $|Q\cap Q_l|=0$, it follows $\int_Q|a_l|=\int_{Q\cap Q_l}|a_l|=0$. Therefore,
 \begin{equation*}\label{ }
Ma_l(x)=\sup_{x\in Q}\frac{1}{|Q|}\int_Q|a_l(y)|dy=\sup_{x\in Q,l(Q)> |x-x_0|/2}\frac{1}{|Q|}\int_{Q\cap Q_l}|a_l(y)|dy\leq \frac{C\|a_l\|_{L^1}}{|x-x_0|^n}
\end{equation*}
for $x\in (2n^{1/2}Q_l)^c$, that is that $M$ satisfies (1.5). While (1.6) follows from the $L^s$ boundedness of $M$ with $1<s\leq \infty$.

Thus, Theorem 5.1 has been proved.

\begin{theorem}\label{th_4.3}
      Let $0<p < \infty ,   0< s \leq \infty$ and $  w(x) \in A_{\infty}$. If $f\in ML^p_w$, then there exist a sequence $\{a_l\}$ of $(p,s,w)$-blocks
 and a sequence $\{\lambda_l\}$ of real numbers such that
 \begin{equation*}
f(x)=\sum_{l}\lambda_l a_l
\end{equation*}
$w$-a.e. and  in $L^p_w$, and

  \begin{equation*}
\left(\sum_{l}|\lambda_l|^p\right)^{1/p}
\leq
C \|f\|_{ML^p_w}.
 \end{equation*}
 In particular, when $0<p\leq 1,$ then $ f\in B^{p,s}_w$, consequently, $ML^p_w \subset B^{p,s}_{w}$.
\end{theorem}

\par  \proof  Let $f\in ML^p_w$ and define for $k=0,\pm 1,\pm 2,\cdots,$
  $$
 E_k=\{x:Mf(x)>2^k\}.
 $$
 Clearly, $E_k$ is open. The Whitney decomposition theorem (see Stein  \cite{St61970}) provides us with closed dyadic cubes  $ Q_k^j$ with the following properties:

 \par (a) $E_k= \bigcup_{j=1}^{\infty}Q_k^j,k=0,\pm 1,\pm 2,\cdots.$

 \par (b) The interiors of the cubes $Q_k^i$ and $Q_k^j$ are disjoint whenever $i\neq j.$

 \par (c) If $l>k,$ then for each $j$ there is an $i$ such that $Q_l^j\subseteq Q_k^i.$

 \par For each integer $k$, let
 $$
 g_k(x) =f(x) \chi _{{\bf R}^n \backslash E_k}(x),b_k^i(x)=\chi_{Q_k^i}(x)f(x)~for ~i=1,2,\cdots.
 $$
 Then
 $$
 f(x)=g_k(x)+\sum_{i=1}^{\infty}b_k^i(x)
 $$
 a.e. (also $w$-a.e. by Lemma 3.1(i)) and in $L^p_w$ since $|g_k(x)|+\sum_{i=1}^{\infty}|b_k^i(x)|\leq |f(x)|$ a.e. (also $w$-a.e.), which is in  $L^p_w$.

\par  We have also
 $g_k(x)\rightarrow 0$  everywhere as $k\rightarrow -\infty$ since  $|g_k(x)|\leq 2^k \rightarrow 0$ as $k\rightarrow -\infty$  and in  $L^p_w$ since $|g_k(x)|\leq |f(x)|$ a.e. (also $w$-a.e.).
  Also
   \begin{equation}
   g_k(x)\rightarrow f(x),w-a.e,  as~ k\rightarrow \infty,
    \end{equation}
    since  $f(x)-g_k(x)$ lives in the set $\{x:Mf(x)>2^k\} $ which decreases to a set of measure 0 as  $k\rightarrow \infty$
  and in  $L^p_w$ since $|f(x)-g_k(x)|\leq 2|f(x)|$ a.e. (also $w$-a.e.).  In fact, we see for given $\varepsilon>0$ that
  \begin{equation}
  \{x:|f(x)-g_k(x)|\geq \frac{\varepsilon}{2} \}\subset \{x:|f(x)-g_k(x)|\neq 0\}\subset E_k.
  \end{equation}
  From (5.6), we can  prove that
\begin{equation}
  \{x:\sup_{j\geq k}|f(x)-g_j(x)|\geq \varepsilon \}\subset
  \bigcup_{j\geq k} \{x:|f(x)-g_j(x)|\geq \frac{\varepsilon}{2} \}\subset E_k,
 \end{equation}
  in fact, take
  $$
  x\in \{x:\sup_{j\geq k}|f(x)-g_j(x)|\geq \varepsilon \},
  $$
  it follows
  $$
 \sup_{j\geq k}|f(x)-g_j(x)|\geq \varepsilon ,
  $$
  then, there is $j_0\geq k$ such that
   $$
 |f(x)-g_{j_0}(x)|\geq\frac{ \varepsilon}{2} ,
  $$
  by (5.6), it follows
  $$
  x\in \{x:|f(x)-g_{j_0}(x)|\geq\frac{ \varepsilon}{2}\}\subset
  \bigcup_{j\geq k} \{x:|f(x)-g_j(x)|\geq \frac{\varepsilon}{2} \}\subset E_k,
  $$
  then, (5.7) holds.   At the same time, we see that
    \begin{eqnarray*}
\infty
> \int_{{\bf R}^n}(Mf)^pw
\geq \int_{\{x:Mf(x)>2^k \}}(Mf)^pw
\geq \int_{\{x:Mf(x)>2^k \}}2^{kp}w
=2^{kp}w(E_k),
 \end{eqnarray*}
it follows
    $$
  w(E_k)\leq \frac{1}{2^{kp}}\int_{{\bf R}^n}(Mf)^pw
  \rightarrow 0 ~{\rm as}~ k\rightarrow \infty.
  $$
    Then, by (5.7),
  $$
  \lim_{k\rightarrow\infty}
  w(\{x:\sup_{j\geq k}|f(x)-g_j(x)|\geq \varepsilon \}) \leq \lim_{k\rightarrow\infty} w(E_k) =0.
  $$
    (5.5)  follows.
   Thus,
  $$
 f(x)=\sum_{k=-\infty}^{\infty}(g_{k+1}(x)-g_{k}(x))
 $$
  $w$-a.e. and in $L^p_w.$
   Thus,
  \begin{eqnarray*}
 f(x)
 &=&
 \sum_{k=-\infty}^{\infty}\left(\sum_i b^i_{k}(x)-\sum_j b^j_{k+1}(x)\right)
\\
 &=&
 \sum_{k=-\infty}^{\infty}\sum_i \left( b^i_{k}(x)-\sum_{j:Q_{k+1}^j\subseteq Q_{k}^i} b^j_{k+1}(x)\right)
 \\
 & &~~~~~~~~~~~~~~~{\rm since~ the ~above ~property ~(c)~ of}~  Q_k^j
 \\
 &=&
 \sum_{k=-\infty}^{\infty}\sum_i \beta^i_{k}(x),
 \end{eqnarray*}
 where
 $$
 \beta^i_{k}(x)= b^i_{k}(x)-\sum_{j:Q_{k+1}^j\subseteq Q_{k}^i} b^j_{k+1}(x).
 $$
We see that supp$\beta^i_{k}(x) \subseteq Q_{k}^i$ and $\beta^i_{k}(x)=g_{k+1}(x)-g_{k}(x)$ if $x\in Q_{k}^i$, so $|\beta^i_{k}(x)|\leq 3\times 2^k.$
Let
$$a_k^i(x)=(w(Q_{k}^i))^{-1/p}(3\times 2^k)^{-1}\beta^i_{k}(x).
$$
Then $\|a_k^i\|_{L^s}\leq |Q_{k}^i|^{1/s}(w(Q_{k}^i))^{-1/p}$, so $a_k^i$ is a $(p,s,w)$-block.
$$
f(x)=\sum_{k=-\infty}^{\infty}\sum_{i}^{}\lambda_k^i a_k^i ~{\rm with} ~\lambda_k^i=3\times 2^k(w(Q_{k}^i))^{1/p},
$$
and
  \begin{eqnarray*}
\sum_{k}^{}\sum_{i}^{}|\lambda_k^i|^p
&=&
\sum_{k}^{}\sum_{i}^{}3^p\times 2^{kp}w(Q_{k}^i)
\\
&=&
3^p\sum_{k}^{} 2^{kp}w(\{x:Mf(x)>2^k\})
\\
&\leq &
C \int_{0}^{\infty} \lambda^{p-1}w(\{x:Mf(x)>\lambda\})d\lambda
\\
&= &
C \|f\|_{ML_w^p}^p.
 \end{eqnarray*}
Thus, we finish the proof of Theorem 5.2.

\par \proof [Proof of Theorem 1.9]
 Noticing $  RH_{r}\subset A_{\infty}$ for $1<r<\infty$,   (1.4) follows from   Theorem 5.1, Theorem 5.2 and  Theorem 1.8. While $ML^1_1 \cap L^1 = \{0\}$, see \cite{Grafakos}, it is clearly not equal to $ H^1$,
  i.e. (1.4) does  not hold when $ q= p =1$ and $w=1$ which is in $A _{1,r}   $ for $1<r<\infty$. Thus, we have proved Theorem 1.9.

\par \section
 { Molecular characterization  of $BH^{p,s}_{w}$}

Coifman \cite{Coif}, Coifman and Weiss \cite{CW}, and Taibleson and
Weiss  \cite{TW2} proved that Hardy spaces $H^p$ with$0<p\leq 1$ can  also be characterized   in terms of the molecular.
Lee and Lin \cite{LL} proved that the weighted  Hardy spaces $H_w^p$ with$0<p\leq 1$ are  characterized   in terms of the weighted molecular.
The molecular characterization of Hardy spaces  provides an effective method to prove boundedness of operators on Hardy spaces.

\par In this section, we give a molecular characterization of $BH^{p,s}_{w}$ that will be used latter.

 \par In the rest of this article, we will denote by $Q^{x_0}_{l} $ the cube centered at $x_0$ with side length $2l$ and denote $Q^0_l $ simply by $Q_l$.

\begin{definition}\label{def 5.1}
  Let $0<p <\infty, w\in A _{q,r}     $ with $ 0<q <p  $  and  $ 1<r <\infty  $,
     and $\max\{ rp/(r-1), p/q\}< s\leq \infty$.
  Set $0<\varepsilon <1 -q/p,
a =1-q/p- \varepsilon,$ and $ b=1-1/s- \varepsilon.$
  A function $M(x)\in L^s $
 is said to be a  $ (p, s,q,w, \varepsilon )$-molecule (centered at  $x_0$), if
\begin{equation}
\|M\|^{a/b}_{L^{s}   }
\left \|
M(x)\left( |Q^{x_0}_{|x-x_0|}|^{-1/s}w(Q^{x_0}_{|x-x_0|})^{1/p}\right)^{b/(b-a)}
\right \|^{1-a/b}_{L^{s}   }
\equiv \Re (M)<\infty.
\end{equation}

\end{definition}

\begin{theorem}\label{th_5.3}
  Let $ p, s, q, w, \varepsilon$ be as in Definition 6.1. Then,
   every  $ (p, s, q, w, \varepsilon )$-molecular $M(x)$  centered
at any point is in $BH^{p,s}_{w}$ and  $\|M \|^{
}_{BH^{p,s}_{w}}\leq C \Re(M)$,
 where the constant $C$ is independent of $M$.
\end{theorem}

\par \proof   Without loss of generality, we can assume
$\Re(M)=1$. In fact, assume $\|M \|^{ }_{BH^{p,s} _{w}}\leq C  $ holds whenever $\Re(M)=1$.
Then, for general $M,$ let $M'=M/\Re(M)$. We have  $\Re(M')=1$ and hence
$\|\Re(M)M' \|^{ }_{BH^{p,s}_{w}}\leq \Re(M)\| M' \|^{}_{BH^{p,s}_{w}}\leq C \Re(M)$.
\par
Let $M$ be a $ (p, s, q, w, \varepsilon )$-molecular centered at
$ x_0$
 satisfying $\Re(M)=1$. Define $Q^{x_0}_l$ by setting
 \begin{equation}\label{5.3}
 \|M\|^{ }_{L^{s}   }= |Q^{x_0}_l|^{1/s}w(Q^{x_0}_l)^{-1/p}.
  \end{equation}
 Let
$2^{k_0-1}<l\leq2^{k_0}$,
 and consider the sets
 $$ E^{x_0}_0=Q^{x_0}_{2^{k_0}}, ~~~~
  E^{x_0}_k=Q^{x_0}_{2^{k_0+k}} \setminus  Q^{x_0}_{2^{k_0+k-1}}  ,~~ for ~~ k=1, 2, \cdots
 ,$$
 Set
 $$M_k=M{\chi _{E^{x_0}_k}},~~ ~~ k=0,1, 2, \cdots.$$
 Let $a,b$ be as in Definition 6.1, from $1<p/q<s$, we see that
 $a>0,b>0, b-a>0.$
By $\Re (M)=1$, we have from (6.1) and (6.2)
 that
 \begin{eqnarray}\label{5.4}
\left \|
M(x)\left( |Q^{x_0}_{|x-x_0|}|^{-1/s}w(Q^{x_0}_{|x-x_0|})^{1/p}\right)^{b/(b-a)}
\right \|_{L^{s}   }
\end{eqnarray}
\begin{eqnarray*}
~~=
|Q^{x_0}_{l}|^{-\frac{1}{s}\frac{a}{b-a}}w(Q^{x_0}_{l})^{\frac{1}{p}\frac{a}{b-a}}
 \leq
C \left( |Q^{x_0}_{2^{k_0}}|^{-\frac{1}{s}}w(Q^{x_0}_{2^{k_0}})^{\frac{1}{p}}\right)^{\frac{a}{b-a}},
 \end{eqnarray*}
noticing $2^{k_0-1}<l\leq2^{k_0}$.

By   the right inequality of (3.1) since $w\in RH_{r}$, and noticing
  $rp/(r-1)<s\leq \infty$, we have
$$
\frac{w(Q^{x_0}_{R_1})}{w(Q^{x_0}_{R_2})}
\leq \left(\frac{|Q^{x_0}_{R_1}|}{|Q^{x_0}_{R_2}|}\right)^{(r-1)/r}
\leq \left(\frac{|Q^{x_0}_{R_1}|}{|Q^{x_0}_{R_2}|}\right)^{p/s}
$$
for $R_1\leq R_2$, it follows
 \begin{equation}\label{5.2}
|Q^{x_0}_{R_1}|^{-1/s}w(Q^{x_0}_{R_1})^{1/p}
\leq
|Q^{x_0}_{R_2}|^{-1/s}w(Q^{x_0}_{R_2})^{1/p}.
 \end{equation}
 Let $x\in E^{x_0}_k $,    we see that $|x-x_0|\geq 2^{k+k_0-1}$.  For   $k=1,2,\cdots,$  we have
\begin{eqnarray*}
& &\int _{E_k} M_k^s(x)dx
\\
 &\leq &
 \int _{E_k} M_k^s(x)
 \left (
 \frac
 {|Q^{x_0}_{|x-x_0|}|^{-1}
 w(Q^{x_0}_{|x-x_0|})^{s/p}}
 {|Q^{x_0}_{2^{k+k_0-1}}|^{-1}
 w(Q^{x_0}_{2^{k+k_0-1}})^{s/p}}
  \right)
  ^{b/(b-a)} dx
  \\
  &&
  ~~~~~~~~ ({\rm since}~ (6.4))
\\
&\leq &
C
 \int _{E_k} M_k^s(x)
 \left (
 \frac
 {|Q^{x_0}_{|x-x_0|}|^{-1}
 w(Q^{x_0}_{|x-x_0|})^{s/p}}
 {|Q^{x_0}_{2^{k+k_0}}|^{-1}
 w(Q^{x_0}_{2^{k+k_0}})^{s/p}}
  \right)
  ^{b/(b-a)} dx
  \\
  &&  ~~~~~~~~~~~~~~~~~~~~~~~~~~~( {\rm by~ }(1.1) )
    \\
&\leq &
C
\left(
|Q^{x_0}_{2^{k+k_0}}|^{1/s}
 w(Q^{x_0}_{2^{k+k_0}})^{-1/p}
 \right) ^{bs/(b-a)}
 \left(
 |Q^{x_0}_{2^{k_0}}|^{1/s}w(Q^{x_0}_{2^{k_0}})^{-1/p}
 \right)^{-as/(b-a)}
 \\
 & &
 ~~( {\rm by}~(6.3))
 \\
&=&
C
\left(
|Q^{x_0}_{2^{k+k_0}}|^{1/s}
 w(Q^{x_0}_{2^{k+k_0}})^{-1/p}
 \right) ^{bs/(b-a)}
 \left(
 |Q^{x_0}_{2^{k+k_0}}|^{1/s}w(Q^{x_0}_{2^{k+k_0}})^{-1/p}
 \right)^{-as/(b-a)}
 \\
 & &
 \times
 \left(
  \frac
 {|Q^{x_0}_{2^{k+k_0}}|^{1/s}w(Q^{x_0}_{2^{k+k_0}})^{-1/p}}
 {|Q^{x_0}_{2^{k_0}}|^{1/s}w(Q^{x_0}_{2^{k_0}})^{-1/p}}
 \right)^{as/(b-a)}
 \\
&\leq&
C
2^{-kn[(r-1)/r-p/s]as/[(b-a)p]}
\left(
|Q^{x_0}_{2^{k+k_0}}|^{1/s}
 w(Q^{x_0}_{2^{k+k_0}})^{-1/p}
 \right) ^{s}
,
\\
&&
 ~~~~~~( {\rm by~the ~right ~inequality~ of }~(3.1)).
\end{eqnarray*}

\par For $s=\infty$ and  $k=1,2,\cdots,$ we see that
\begin{eqnarray*}
\|M_k\|_{L^{\infty}}
 &=&
 \|
 M_k(x)
  w(Q^{x_0}_{|x-x_0|})^{b/[p(b-a)]}
   w(Q^{x_0}_{|x-x_0|})^{-b/[p(b-a)]}
 \|_{L^{\infty}}
\\
& \leq &
  w(Q^{x_0}_{2^{k+k_0-1}})^{-b/[p(b-a)]}
\|
 M_k(x)
  w(Q^{x_0}_{|x-x_0|})^{b/[p(b-a)]}
  \|_{L^{\infty}}
 \\
 &&
 ~~~~( {\rm by }~(6.4))
\\
& \leq &
  w(Q^{x_0}_{2^{k+k_0}})^{-b/[p(b-a)]}
 w(Q^{x_0}_{2^{k_0}})^{a/[p(b-a)]}
 \\
 &&
  ~~~~~~~( {\rm  by ~(1.1)~ and  } ~(6.3))
  \\
& \leq &
C2^{-n(r-1)ak/[rp(b-a)]}
  w(Q^{x_0}_{2^{k+k_0}})^{-1/p}
   \\
 &&~~~~~~~~~~~~( {\rm by~the ~right ~inequality~ of  }~ (3.1)).
\end{eqnarray*}
Noticing $2^{k_0-1}<l\leq2^{k_0}$,  by  (6.2),  (6.4),  and  (1.1),
 we have
\begin{eqnarray*}
\|M_0\|^{ }_{L^{s}   }
&\leq& |Q^{x_0}_l|^{1/s}w(Q^{x_0}_l)^{-1/p}
\\
&\leq & |Q^{x_0}_{2^{k_0-1}}|^{1/s}w(Q^{x_0}_{2^{k_0-1}})^{-1/p}
\leq C |Q^{x_0}_{2^{k_0}}|^{1/s}w(Q^{x_0}_{2^{k_0}})^{-1/p}
.
\end{eqnarray*}
 Hence, for
 $k=0, 1,2,\cdots,$ we have

$$\|M_k\|^{ }_{L^{s}   } \leq C 2^{-\delta k}
 |Q^{x_0}_{2^{k+k_0}}|^{1/s}
 w(Q^{x_0}_{2^{k+k_0}})^{-1/p}
$$
 where $\delta=n[\frac{r-1}{r}-\frac{p}{s}]\times \frac{a}{(b-a)p}>0 $ that follows from $rp/(r-1)<s\leq \infty$, $C$ is a constant, let
$$
a_k(x)=\frac{1}{C}2^{\delta k}M_k(x), ~~k=0,1,2,\cdots,
$$
 then
 \begin{equation}\label{6.6}
M(x)=\sum \limits_{k=0}^{\infty}M_k(x)
=\sum \limits_{k=0}^{\infty}C2^{-\delta k}a_k(x), {\rm ~for~ all~} x\in {\bf R}^n,
  \end{equation}
 and each $a_k$ is a $ (p, s, w )$-block centered at  $x_ 0$ with
${\rm supp}a_k \subset Q^{x_0}_{2^{k+k_0}}
,$ and $\sum
\limits_{k=0}^{\infty}C2^{-\delta \bar{p}k} =C<\infty$.
 (6.5) holds for all $x\in {\bf R}^n$, and also in $L^p_w$ by the dominated convergence theorem.
Thus, $M\in BH^{p,s}_{w}$ and
$\|M \|^{ }_{BH^{p,s}_{w}}\leq C $.

\par  Thus, we complete the proof of the theorem.

\par On the one hand, we have

\begin{proposition}\label{rem 5.2}
  Every $(p,s,w)$-block $h$ is a $(p,s,q,w,\varepsilon)$-molecule for $p,s, q , w,\varepsilon $ in Definition 6.1.
\end{proposition}

\par \proof

Let $h$ be a $(p,s,w)$-block with supp $h \subset Q^{x_0}_R$ and $ \|h\|_{L^s}\leq | Q^{x_0}_R|^{1/s}w( Q^{x_0}_R)^{-1/p}$.
Since $|x-x_0|\leq R$ for $x\in Q^{x_0}_R$, then, by (6.4),
 we have
\begin{eqnarray*}
& &
\left \|
h(x)\left( |Q^{x_0}_{|x-x_0|}|^{-1/s}w(Q^{x_0}_{|x-x_0|})^{1/p}\right)^{b/(b-a)}
\right \|_{L^{s}   }
\\
&\leq &
\left( |Q^{x_0}_{R}|^{-1/s}w(Q^{x_0}_{R})^{1/p}\right)^{b/(b-a)}
\left \|h\right \|_{L^{s}   },
\end{eqnarray*}
it follows $  \Re (h)  \leq   |Q^{x_0}_{R}|^{-1/s}w(Q^{x_0}_{R})^{1/p}\|h\|^{}_{L^{s}}   \leq C$, by the definition of $h$. Thus, we finish the proof of the proposition.

\par \section
 { Estimates of  operators for $0<p<\infty$
 }

In this section, we prove   Theorem 1.10, 1.12, 2.2, 2.4 and 2.9.

\par \subsection
 { Proof of Theorem 1.10
  }

\par  \proof[Proof of Theorem 1.10] By Lemma 3.8, we have  that $\|T h\|_{L^{p}_{w}}\leq C$ for every $(p,s,w)$-block $h$.
  Replacing $a_i$ with $Ta_i$, repeating the proof process of conclusion (ii) of Proposition 3.3, and using Proposition 3.6,  we get  Theorem 1.10(i). By (3.3), we get  Theorem 1.10(ii).

\par \subsection
 { Proof of Theorem 1.12
   }

\par
\proof[Proof of Theorem 1.12] Let $ f= \sum_{j=1}^{\infty}\lambda_ja_j \in BH^{p,s}_{w}$ where each $a_j$ is a $(p,s,w)$-block and $\sum_{j=1}^{\infty}|\lambda_j|^{\bar{p}}<\infty$.
We claim that there is a constant $C$ such that, for $w\in A_{q,r}$ with $0<q<p$ and $1<r<\infty$, for every  $(p,s,w)$-block $a_j$,
\begin{equation}\label{7.5}
\|T a_j\|_{BH^{p,s}_{w}}\leq C.
\end{equation}
Once (7.1) is proved, for a linear operator $T$,  using  $\sum_{j=1}^{\infty}|\lambda_j|^{\bar{p}}<\infty$,   the subadditivity and the completeness  of $BH^{p,s}_{w}$ (Proposition 3.5 and Proposition 3.6),   we have by a standard argument that $T$ has an unique bounded extension (still denoted by $T$) from $BH^{p,s}_{w}$     to itself satisfying (1.8)  in $BH^{p,s}_{w}$ and $w$-a.e.. By (1.8), (3.7) and (7.1),  we get (1.11). The part (i) of the theorem is proved.

\par For a sublinear operator $T$  which satisfies (1.7), noticing  $w\in  A_{q,r}\cap P $, we see from Lemma 3.1 that  (1.7) holds a.e.. Then,
using the monotonicity (Proposition 3.4), the subadditivity (3.7)  (Proposition 3.5) of $BH^{p,s}_{w}$ , and (7.1),   we get (1.11). The part (ii) of the theorem is proved.

\par Next, Let us prove (7.1). By Theorem 6.2,
it  suffices to check that $Th$ is a  $ (p, s, q, w,\varepsilon
)$-molecular  for every $(  p,s,w)-$block $h$ and $\Re (Th)\leq C $ for a constant  $C$  independent of $h$. To do this,
 we suppose   supp $h\subset Q^{x_0}_{2^{m_0}} $ and $\|h\|_{L^{s}}\leq C |Q^{x_0}_{2^{m_0}}|^{1/s}w(Q^{x_0}_{2^{m_0}})^{-1/p}$, where $m_0 \in {\bf R} $. Let $2^{l_0}=n^{1/2}, k_0=m_0+l_0, $ then, $ n^{1/2}Q^{x_0}_{2^{m_0}}=Q^{x_0}_{2^{k_0}} .$
Let  $\varepsilon, a$ and $b$ as in Definition 6.1.   We have   $a>0,b>0$ and $b-a> 0$.

\par  To prove $\Re (Th)\leq C $, we write, for $1< s<\infty$,
\begin{eqnarray*}
J=: \left\|Th(x)\left( |Q^{x_0}_{|x-x_0|}|^{-1/s}w(Q^{x_0}_{|x-x_0|})^{1/p}\right)^{b/(b-a)}\right\|^{s}_{L^{s}}  =:  J_1+J_2,
\end{eqnarray*}
where
\begin{eqnarray*}
J_1 =
 \int _{Q^{x_0}_{2^{k_0+1}}}
|Th(x)|^s \left( |Q^{x_0}_{|x-x_0|}|^{-1/s}w(Q^{x_0}_{|x-x_0|})^{1/p}\right)^{bs/(b-a)} dx,
\end{eqnarray*}
\begin{eqnarray*}
J_2 =
 \int _{(Q^{x_0}_{2^{k_0+1}})^c}
|Th(x)|^s \left( |Q^{x_0}_{|x-x_0|}|^{-1/s}w(Q^{x_0}_{|x-x_0|})^{1/p}\right)^{bs/(b-a)} dx,
\end{eqnarray*}
and,  for $s=\infty,$
\begin{eqnarray*}
I=:\left\|Th(x)\left( w(Q^{x_0}_{|x-x_0|})^{1/p}\right)^{b/(b-a)}\right\|_{L^{\infty}} =: \max\{ I_1,I_2\},
\end{eqnarray*}
where
\begin{eqnarray*}
I_1=
\left\|Th(x)\left( w(Q^{x_0}_{|x-x_0|})^{1/p}\right)^{b/(b-a)}\right\|_{L^{\infty}(Q^{x_0}_{2^{k_0+1}})},
 \end{eqnarray*}
\begin{eqnarray*}
I_2=
\left\|Th(x)\left( w(Q^{x_0}_{|x-x_0|})^{1/p}\right)^{b/(b-a)}\right\|_{L^{\infty}((Q^{x_0}_{2^{k_0+1}})^c)}.
 \end{eqnarray*}
\par We need to estimate $J_i$ and $ I_i,i=1,2$. Noticing $w\in RH_r$ and  $ rp/(r-1)<s\leq \infty$,
by using (6.4) ($R_1=|x-x_0|, R_2=2^{k_0+1}$),   $L^s$-boundedness of $T$, and (1.1), we have
$$
J_1
 \leq  C \left(
|Q^{x_0}_{2^{k_0}}|^{-1/s}w(Q^{x_0}_{2^{k_0}})^{1/p}
\right)^{bs/(b-a)}
 \|h \|^{s}_{L^{s}   },
 $$
and
$$
I_1
 \leq  C \left(
w(Q^{x_0}_{2^{k_0}})^{1/p}
\right)^{b/(b-a)}
 \|h \|^{}_{L^{\infty}   }.
 $$
To estimate $J_2$ and $I_2$,  let  $x\in (Q^{x_0}_{2^{k_0+1}})^c$, we see
 $ Q^{x_0}_{2^{k_0-1}} \subset Q^{x_0}_{|x-x_0|}$.    Since $w\in D_q$, we have by (1.1) that
\begin{equation}\label{7.9}
w( Q^{x_0}_{|x-x_0|})
\leq 2^{nq}(2^{-k_0}|x-x_0|)^{nq}w(Q^{x_0}_{2^{k_0}}).
\end{equation}
\par For $J_2$, notice that  $ p/q<s<\infty$,
 we have
\begin{eqnarray*}
& J_2&
\\
& \leq &
C
\int _{(Q^{x_0}_{2^{k_0+1}})^c}
|Th(x)|^s \left( |x-x_0|^{-n/s}
(2^{-k_0}|x-x_0|)^{nq/p}
w(Q^{x_0}_{2^{k_0}})^{1/p}\right)^{bs/(b-a)} dx
\\
& &
~~~~~~~~~~~~~~~~~~~~~~~~({\rm by}~(7.2))
\\
& =&
C
\left(
2^{-k_0nq}
w(Q^{x_0}_{2^{k_0}})\right)^{bs/[p(b-a)]}
\int _{(Q^{x_0}_{2^{k_0+1}})^c}
|Th(x)|^s  |x-x_0|^{bns}
 dx
\\
& &
~~~~~~~~~~~~~~~~~~~~~~~~({\rm noticing} ~b-a=q/p-1/s)
\\
& \leq &
C
\left(
|Q^{x_0}_{2^{k_0}}|^{-q}
w(Q^{x_0}_{2^{k_0}})\right)^{bs/[p(b-a)]}
\|h\|_{L^1}^s
\int _{(Q^{x_0}_{2^{k_0+1}})^c}
  |x-x_0|^{(b-1)ns}
 dx
\\& &~~~~~~~~~~~~~~~~~~~
~~~~~({\rm by}~(1.5))
\\
& \leq &
C
\left(
|Q^{x_0}_{2^{k_0}}|^{-q}
w(Q^{x_0}_{2^{k_0}})\right)^{bs/[p(b-a)]}
\|h\|_{L^s}^s |Q^{x_0}_{2^{k_0}}|^{(1-1/s)s}
|Q^{x_0}_{2^{k_0}}|^{s(b-1+1/s)}
 \\
& &
~~~~~~~~~~~~~~~({\rm by ~H\ddot{o}lder~inequality~and ~
noticing }~b-1+1/s=-\varepsilon <0)
\\
& = & C \left(
|Q^{x_0}_{2^{k_0}}|^{-1/s}w(Q^{x_0}_{2^{k_0}})^{1/p}
\right)^{bs/(b-a)}
 \|h \|^{s}_{L^{s}   }
 \\& &~~~~~~~~~~~~~
 ~~~~~~~~~~~({\rm noticing} ~b-a=q/p-1/s).
\end{eqnarray*}
For  $I_2,$ we have for $x\in (Q^{x_0}_{2^{k_0+1}})^c$ that
\begin{eqnarray*}
& & |Th(x)|\left( w(Q^{x_0}_{|x-x_0|})^{1/p}\right)^{b/(b-a)}
\\
&\leq &
C\frac{\|h\|_{L^1}}{|x-x_0|^n}
\left(\left((2^{-k_0}|x-x_0|)^{nq}w(Q^{x_0}_{2^{k_0-1}})\right)^{1/p}\right)^{b/(b-a)}
\\
& &~~~~~~~~~~~~
~~~~~~~~~~~~ ({\rm by}~(1.5)~{\rm and}~(7.2))
\\ &= & C
\left(\left(2^{-k_0nq}w(Q^{x_0}_{2^{k_0-1}})\right)^{1/p}\right)^{b/(b-a)}\frac{\|h\|_{L^1}}{|x-x_0|^{n(1-b)}}
\\&&~~~~~~~~~~~~~~~~~~~~~~~~{\rm noticing }~b-a=q/p)
\\
&\leq  &
C
\left(\left(2^{-k_0nq}w(Q^{x_0}_{2^{k_0-1}})\right)^{1/p}\right)^{b/(b-a)}
2^{-k_0n(1-b)}
\|h\|_{L^1}
\\
&&~~~~~~~~~~~~~~~~~~~~~~~~({\rm since} ~x\in (Q^{x_0}_{2^{k_0+1}})^c ~{\rm and} ~ 1-b=\varepsilon>0)
\\
&\leq  &
C
\left(\left(w(Q^{x_0}_{2^{k_0}})\right)^{1/p}\right)^{b/(b-a)}
\|h\|_{L^{\infty}}
\\&&~~~~~~~~~~~~~~~~~~
~~~~~~({\rm noticing} ~b-a=q/p),
\end{eqnarray*}
thus,
$I_2
\leq
C\left(\left(w(Q^{x_0}_{2^{k_0}})\right)^{1/p}\right)^{b/(b-a)}
\|h\|_{L^{\infty}}.
$

\par
Combining with the above estimate of $J_i$ and $I_i, i=1,2$, by using (1.1), and noticing $b>0$ and $b-a>0$,
  we have  that
\begin{eqnarray*}\label{7.10}
&&\left\|Th(x)\left( |Q^{x_0}_{|x-x_0|}|^{-1/s}w(Q^{x_0}_{|x-x_0|})^{1/p}\right)^{b/(b-a)}\right\|_{L^{s}}
\\
&&\leq
C
\left(
|Q^{x_0}_{2^{m_0}}|^{-1/s}w(Q^{x_0}_{2^{m_0}})^{1/p}
\right)^{b/(b-a)}
 \|h \|^{}_{L^{s}   }
\end{eqnarray*}
 for  $ \max \{ rp/(r-1),q/p\}<s\leq \infty$.
 At the same time, by the boundedness of $T$ on $L^s$,
$\|Th \|^{ }_{L^{s}   } \leq  C \| h  \|^{ }_{L^{s}   } .$
Then,  we have
\begin{eqnarray*}
\Re (Th)\leq C\| h  \|^{ }_{L^{s}   }|Q^{x_0}_{2^{m_0}}|^{-1/s}w(Q^{x_0}_{2^{m_0}})^{1/p}\leq C.
\end{eqnarray*}
Thus, we complete the proof of Theorem 1.12.

\par \subsection
 {Proof of  Theorem 2.1, 2.2, 2.4 and 2.9}

\par \proof[Proof of Theorem 2.1] Theorem 2.1 follows from the following lemma.
\begin{lemma}\label{pro 7.12}  (2.1) implies (1.5).
\end{lemma}

\proof
Let $h$ be a $(p,s,w)$-block.
We suppose that $\supp h\subset Q$ with centre $x_0$. Take $x\in (2n^{1/2}Q)^c$.  For any $y\in Q$,   we see
 that   $|x-x_0|\leq 2|x-y|$. Then, by (2.1), $|Th(x)|\leq C 2^n\frac{\|h\|_{L^1}}{|x-x_0|^n}$
for $x\in (2n^{1/2}Q)^c$, i.e. (1.5) holds. The lemma holds.

\par \proof[Proof of Theorem 2.2]
Clearly, each  $T_{1}$  is  linear, and   satisfy (2.1) (see also  \cite {SW} and \cite {DingL}),  by Lemma 7.1, (1.5) follows.
  It is known that each  $T_{1}$ is bounded on  $L^s$ with $1<s<\infty$, and (1.6) follows.
  Then the theorem follows from   Theorem 1.10(i) and 1.12(i).

\par \proof [Proof of Theorem   2.4]
 Each $T_2$ obviously  satisfies (2.1),   by Lemma 7.1, (1.5) follows, and  the corresponding  constant $C$ in (2.1) and (1.5) does not depend on the associated $\varepsilon, \xi, P$, et al, therefore, each $T_2^*$   satisfies  (1.5).

  It is known that each $T_2^*$ is bounded on $L^s$  for all $1<s<\infty$, (1.6) follows.

\par Each $T_2$  satisfies   (1.5) and is linear. It is also  bounded on $L^s$  for all $1<s<\infty$,  then it satisfies (1.6). By Theorem 1.10(i),   $T_2$ has an unique bounded extension (still denoted by $T_2$) from $BL^{p,s}_{w}$
    to $L^{p}_{w}$ that satisfies (1.8)  $w$-a.e., it follows that the corresponding   maximal operators  $T_2^*$ satisfies (1.7).

\par  Each $T_2^*$ obviously is sublinear.

 By     Theorem 1.10(ii) and Theorem 1.12(ii), then Theorem 2.4 holds.

\par \proof [Proof of Theorem 2.9] We see in the proof of Theorem 5.1 that $M$ satisfies (1.5).
It is known that each  $M$ is bounded on  $L^s$ with $1<s<\infty$,  (1.6) follows.
We also  see in Section 5 that $M$ satisfies (1.7). Then the theorem follows from Theorem 1.10(ii) and Theorem 1.12(ii).

\end{document}